\documentclass[11pt]{article}

\usepackage{color}
\usepackage{latexsym}
\usepackage{amssymb}
\usepackage{amsthm}
\usepackage{amsmath}
\usepackage{graphicx}
\usepackage{enumerate}
\usepackage{hyperref}
\usepackage{dsfont}
\usepackage{pict2e}
\usepackage{epic}
\newtheorem{Theorem}{Theorem}[part]
\newtheorem{Definition}{Definition}[part]
\newtheorem{Proposition}{Proposition}[part]

\newtheorem{Lemma}{Lemma}[part]
\newtheorem{Corollary}{Corollary}[part]
\newtheorem{Remark}{Remark}[part]

\def \Frac{\displaystyle\frac}

\def \trans{^{\scriptscriptstyle{\intercal }}}

\def \N{\mathbb{N}}
\def \R{\mathbb{R}}

\def \E{\mathbb{E}}
\def \F{\mathbb{F}}

\def \P{\mathbb{P}}
\def \Q{\mathbb{Q}}

\def \S{\mathbb{S}}

\def \Ac{{\cal A}}
\def \Bc{{\cal B}}

\def \Dc{{\cal D}}

\def \Fc{{\cal F}}
\def \Gc{{\cal G}}
\def \Hc{{\cal H}}
\def \Ic{{\cal I}}
\def \Kc{{\cal K}}
\def \Lc{{\cal L}}

 \def \Nc{{\cal N}}

\def \Vc{{\cal V}}

\def \Vc{{\cal V}}

\def \ni{\noindent}

\def \eps{\varepsilon}

\def \ep{\hbox{ }\hfill$\Box$}

\def\Dt#1{\Frac{\partial #1}{\partial t}}

\def\reff#1{{\rm(\ref{#1})}}

\def\beqs{\begin{eqnarray*}}
\def\enqs{\end{eqnarray*}}
\def\beq{\begin{eqnarray}}
\def\enq{\end{eqnarray}}

\begin{document}
\title{Optimal Switching  in Finite Horizon \\ under State Constraints}
\author{
Idris KHARROUBI\\
\footnotesize{CEREMADE}\\
\footnotesize{CNRS UMR 7534}\\
\footnotesize{Universit\'e Paris Dauphine}\\
\footnotesize{and CREST}\\
\footnotesize{ \texttt{kharroubi @ ceremade.dauphine.fr}
}\\
}

\date{First version: March 2015\\
This version: January 2016}

\maketitle

\begin{abstract}
We study an optimal switching problem with a state constraint: the controller is  only allowed to choose strategies that keep the controlled diffusion in a closed domain. We prove that the value function associated with this problem is the limit of value functions associated with unconstrained switching problems with penalized coefficients, as the penalization parameter goes to infinity. 
This convergence allows to set a dynamic programming principle for the constrained switching problem. We then prove that the value function is a solution to a system of variational inequalities (SVI for short) in the  constrained viscosity sense. 
We finally prove that uniqueness for our SVI cannot hold and we give a weaker characterization of the value function as the maximal solution to this SVI.  All our results are obtained without any regularity assumption on the constraint domain.
\end{abstract}

\vspace{2mm}


\vspace{2mm}


\section{Introduction}
\setcounter{equation}{0} \setcounter{Assumption}{0}
\setcounter{Theorem}{0} \setcounter{Proposition}{0}
\setcounter{Corollary}{0} \setcounter{Lemma}{0}
\setcounter{Definition}{0} \setcounter{Remark}{0}

Optimal control of multiple switching regimes consists in looking for the value of an optimization problem where the allowed strategies are sequences of interventions. It naturally arises in many applied disciplines where it is not realistic to assume that the involved quantities can be continuously controlled. 
More precisely, the optimal switching problem supposes that the control strategies are sequences $\alpha=(\tau_k,\zeta_k)_{k}$ where the sequence $(\tau_k)_{k}$ represents the intervention times of the controller and $\zeta_k$ corresponds to the level of intervention of the agent at each time $\tau_k$. 
 
Such a class of strategies allows to consider discrete actions for the controller which can be more relevant than continuous time controls. Therefore, the modelization with optimal switching problems has attracted a lot of interest during the last decades
(see e.g. Brennan and Schwarz \cite{BS85} for resource extraction, Dixit \cite{D} for
production facility problems, Carmona and Ludkovski \cite{carlud05} for power plant management or Ly Vath,  Pham and Villeneuve \cite{lpv} for dividend decision problem with reversible technology investment). 

Another specificity  to take into account in the modelization with optimal switching is the limitation of the quantities involved in the control problem. Indeed, in most of management problems the controlled system is subject to a constraint on the possible states that it can take. For example, a solvency condition is usually imposed to the investors of a financial market and the energy producer has to take into account the limited storage capacities.
This leads to impose a state constraint on the controlled diffusion $X$ of the form 
\beqs
X_s & \in & \Dc\quad \mbox{ for all } s,
\enqs
 where  $\Dc$ is a closed set. 
 We therefore need to restrict our control problem to the set $\Ac_{t,x}^{\Dc}$  of strategies that keep the controlled diffusion starting from $(t,x)$ in the constraint domain $\Dc$. 
   Unfortunately, such a constraint leads to strong difficulties due, in particular, to the complicated structure of the set valued function $(t,x)\mapsto\Ac_{t,x}^{\Dc}$. 
To the best of our knowledge, no rigorous study of the optimal switching problem in the constrained case has been done before and our aim is to fill this gap.

In the continuous time control case, H. M. Soner gives in \cite{Son86} a first study of the constrained problem in a deterministic framework where he introduces the notion of constrained viscosity solutions.  To characterize the value function, his approach relies on a continuity argument under an assumption on the boundary of the constraint domain $\partial \Dc$. He then extends this result to the case of piecewise deterministic processes in \cite{Son86b}.
The continuous time stochastic control case is studied  by  M. A. Katsoulakis in  \cite{Kat}. 
His approach also relies on continuity and he imposes regularity conditions on the constraint domain $\Dc$.
In our case, such an approach is not possible since the value function may be discontinuous even for a smooth domain $\Dc$ as shows the counter-example presented in Sub-section \ref{contreexemple}.

Let us also mention the recent approach of D. Goreac \textit{et al.} presented in \cite{GIS13}. They formulate the initial problem as a linear problem which concerns the occupation measures induced by the controlled diffusion processes.  Under convexity assumptions, the authors characterize (see Theorem 11 in \cite{GIS13})  the value function associated to the weak formulation of the continuous time stochastic control problem under state constraints (the weak formulation means that the controller is allowed to choose  the probability space in addition to the control strategy). Unfortunately, such an approach cannot be applied to the optimal switching under state constraints since the the set of values taken by the controls is not convex.


In this work, we present an original approach which 
allows to deal with the lack of regularity of the associated value function.
Moreover, our method does not need any regularity or convexity assumption.
In particular, we only need to assume that the constraint domain $\Dc$ is closed.

 To be more precise, our approach relies on the simple structure of switching controls.  Indeed, they can be seen as random variables taking values in $([0,T]\times\Ic)^\N$ where $\Ic$ is a finite set and $T>0$ is a given constant. From Tychonov theorem, we get the compactness of this space which allows to prove the tightness of a sequence $(\alpha^n)_n$ of switching strategies and hence the convergence in law up to a subsequence. 
Then applying Skorokhod representation theorem, we are able to provide a probability space and a sequence $(\tilde \alpha^n)_n$ that converges almost surely to some $\tilde \alpha$ and such that $\tilde \alpha^n$ is equal in law to $\alpha^n$ for all $n$.

We use this sequential compactness property in the following way. We first introduce a sequence $(v_n)_n$ of unconstrained switching problems with $n$-penalized terminal and running reward coefficients out of the constraint domain $\Dc$. For each penalized switching problems $v_n$, we take $\alpha^n$ as a ${1\over n}$-almost optimal strategy for $v_n$ and we make $\tilde \alpha ^n$ converge to $\tilde \alpha$ as described previously. 
Then we construct a switching strategy $\alpha^*$ which is equal in law to $\tilde \alpha$. To this end we prove stability results for measurability and convergence properties for sequence of diffusion driven by converging Brownian motions. These results that have their own interest are presented separately in the Appendix.

The strong convergence of $\tilde \alpha^n$ to $\tilde \alpha$ allows to prove that $\alpha^*$ is optimal for the switching problem under constraint. 
As a byproduct, we get the convergence of the unconstrained penalized switching problems to the constrained one. Using existing results on classical optimal switching problems, this convergence allows to set a dynamic programming principle for the constrained switching problem.



We then focus on the PDE characterization of the value function. Using the dynamic programming principle proved before, we show that the value function is a constrained viscosity solution to a system of variational inequalities (SVI for short) defined on the constraint domain $\Dc$. We then investigate the uniqueness of a solution to this SVI. The usual approach to get uniqueness of a viscosity solution consists in proving a comparison theorem for the PDE. As a consequence of such a comparison theorem, the unique solution has to be continuous.  Unfortunately, the continuity of the value functions is not true in general as shown by the counter-example given in Sub-section \ref{contreexemple}. Therefore, we cannot hope to state such a uniqueness result for the SVI on $\Dc$. Instead, we characterize our value function as the maximal viscosity solution of the SVI under an additional growth assumption. This maximality property is also obtained from the convergence of the penalized unconstrained problems to the constrained one. 

We end the introduction by the description of the organization of the paper.
 In Section 2 we expose in detail the formulation of the optimal switching problem under state constraints and we provide a simple example to stress the possible lack of regularity for the value function. We then give in Section 3 some examples of application. 
 In Section 4, we provide an approximation of our constrained problem by unconstrained problems with penalized coefficients. We prove the convergence of the penalized problems to the constrained one as the penalization parameter goes to infinity. In Section 5, we state a dynamic programming principle and we prove that the value function is a constrained viscosity solution to a SVI. Finally, in section 6 we focus on uniqueness. Since we cannot prove uniqueness of a solution for the SVI, we characterize the value function as the maximal constrained viscosity solution to the SVI under an additional growth assumption. Some examples where this additional growth condition is satisfied are then given. 



\section{Problem formulation}
\setcounter{equation}{0} \setcounter{Assumption}{0}
\setcounter{Theorem}{0} \setcounter{Proposition}{0}
\setcounter{Corollary}{0} \setcounter{Lemma}{0}
\setcounter{Definition}{0} \setcounter{Remark}{0}
\subsection{Optimal switching under state constraints}

We fix a complete probability space $\big(\Omega,\Gc,\P\big)$ which is endowed with a Brownian motion $W=(W_t)_{t\geq0}$ valued in $\R^d$. We denote by $\F$ the complete and right continuous filtration generated by $W$. We also consider a terminal time given by a constant $T>0$.

\paragraph{Controls.} 

We then define the set  $\Ac_t$ of admissible switching controls at time $t\in[0,T]$ as the set of double sequences $\alpha = (\tau_k,\zeta_k)_{k\geq0}$ where 
\begin{itemize}
\item ${(\tau_k)}_{k\geq0}$ is a nondecreasing sequence of  $\mathbb{F}$-stopping times with $\tau_0=t$  and $\lim_{k\rightarrow\infty}\tau_k>T$,
\item  $\zeta_k$ is an
$\Fc_{\tau_k}$-measurable random variables valued in the set $\Ic$ defined by $\Ic=\{1,\ldots,m\}$.
\end{itemize}
With a
strategy $\alpha=(\tau_k,\zeta_k)_{k\geq0}\in\Ac_t$ we associate the
process $(\alpha_s)_{s\geq t}$ defined by
\beqs
\alpha_s & =  & 
\sum_{k\geq0}\zeta_k\mathds{1}_{[\tau_k,\tau_{k+1})}(s)\;,\quad s\geq t\;.
\enqs


\paragraph{Controlled diffusion.} We are given two functions $\mu:~\R^d\times\Ic\rightarrow\R^d$ and $\sigma:~\R^d\times\Ic\rightarrow\R^{d\times d}$. We make the following assumption.

\vspace{2mm}

\ni \textbf{(H1)} There exists a constant $L$ such that
\beqs
|\mu(x,i)-\mu(x',i)|+|\sigma(x,i)-\sigma(x',i)| & \leq & L |x-x'|\;,
\enqs
for all $(x,x',i)\in\R^{d}\times\R^{d}\times\Ic$.

\vspace{2mm}

 For $(t,x)\in[0,T]\times\R^d$ and $\alpha\in\Ac_t$ we consider the controlled diffusion $X^{t,x,\alpha}$ defined by the following SDE
\beq\label{EDS}
X^{t,x,\alpha}_s & = & x+\int_t^s\mu\big(X^{t,x,\alpha}_r,\alpha_r\big)dr+\int_t^s\sigma\big(X^{t,x,\alpha}_r,\alpha_r\big)dW_r\;,\quad  s\geq t\;.
\enq 
Under \textbf{(H1)}, we have existence and uniqueness of an $\F$-adapted solution $X^{t,x,\alpha}$ to \reff{EDS}  for any initial condition $(t,x)\in[0,T]\times\R^d$ and any switching control $\alpha\in\Ac_t$. 

We also have the following classical estimate (see e.g. Corollary 12, Section 5, Chapter 2 in \cite{Kryl80}): for any $q\geq1$ there exists a constant $C_q$ such that 
\beq\label{XS2}
\sup_{\alpha\in\Ac_t}\E\left[\sup_{s\in[t,T]}\big|X^{t,x,\alpha}_s\big|^q\right] & \leq & C_q\big(1+|x|^q\big)
\enq
for all $(t,x)\in[0,T]\times\R^d$.


\paragraph{Expected Payoff.} We consider terminal and running reward functions $g:~\R^d\times\Ic\rightarrow\R$ and $f:~\R^d\times\Ic\rightarrow\R$ and a cost function $c:~\R^d\times\Ic\times\Ic\rightarrow\R$ on which we impose the following assumption.

\vspace{2mm}

\noindent \textbf{(H2)} 

\begin{enumerate}[(i)]
\item The function $f$, $g$ and $c$ are locally Lipschitz: for any $R>0$ there exists a constant $L_R$ such that
\beqs
|g(x,i)-g(x',i)|+|f(x,i)-f(x',i)| +|c(x,i,j)-c(x',i,j)|& \leq & L_R |x-x'|\,,\qquad
\enqs
for all $i,j\in\Ic$ and  $x,x'\in\R^{d}$ such that $|x|\leq R$ and $|x'|\leq R$.

\item  There exists a constant $C$ and an integer $q$ such that 
\beqs
|g(x,i)|+|f(x,i)| +|c(x,i,j)|& \leq & C\big(1+ |x|^q)\,,\qquad
\enqs
for all $x\in\R^{d}$ and $i,j\in\Ic$.

\item  There exists a constant $\bar c$ $>$ $0$, such that 
\beqs\label{majcout}
c(x,i,j) & \geq & \bar c\,,
\enqs
for all $x\in\R^d$ and $i,j\in\Ic$.

\end{enumerate}

\vspace{2mm}
\ni We then define the functional pay-off $J$ up to time $T$ by
\beqs
J\big(t,x,\alpha\big) & = & \E\Big[g\big(X^{t,x,\alpha}_T,\alpha_T\big)+\int_t^Tf\big(X^{t,x,\alpha}_s,\alpha_s\big)ds-\sum_{k\geq 1}c\big(X^{t,x,\alpha}_{\tau_k},\zeta_{k-1},\zeta_{k}\big)\mathds{1}_{\tau_k\leq T}\Big]
\enqs 
for all $(t,x)\in [0,T]\times\R^d$ and $\alpha\in \Ac_t$.

Under \textbf{(H1)} and \textbf{(H2)} we get from \reff{XS2} that $J\big(t,x,\alpha\big)$ is well defined
for any initial condition $(t,x)\in [0,T]\times\R^d$ and any control $\alpha\in \Ac_t$.


\paragraph{State constraint.} Let $\Dc$ be a
nonempty closed subset of $\R^{d}$. For $(t,x,i)\in[0,T]\times\Dc\times \Ic$ we denote by $\Ac_{t,x,i}^{\Dc}$ the set of strategies $\alpha\in  \Ac_t$ such that $\zeta_0=i$ and
\beqs
\P\Big( X^{t,x,\alpha}_s\in \Dc \mbox{ for all }s\in[t,T]  \Big) & = & 1\;.
\enqs

\vspace{2mm}

\paragraph{Value function.}  We then define the value function $v$ associated with the switching problem under state constraints by
\beq\label{opt-switch-const}
v(t,x,i) & = & \sup_{\alpha\in \Ac_{t,x,i}^{\Dc}}J\big(t,x,\alpha\big) 
\enq
for all $(t,x,i)\in[0,T]\times\Dc\times\Ic$, with the convention $v(t,x,i)=-\infty$ if $\Ac_{t,x,i}^{\Dc}=\emptyset$.
Our aim is to give an analytic characterization of the function $v$.
\subsection{Lack of smoothness for the value function}\label{contreexemple}

In general control theory, we expect to get a continuous value function  as we assume that the parameters are continuous.
In the framework of optimal switching under constrains such a property fails to be true. Indeed, the following simple example provides a discontinuous value function.


Fix $d=2$ and consider the case where $\Dc$  is the smooth domain $\R\times\R_+$. Take $\Ic=\{1,2\}$ and define the diffusion coefficients $\mu$ and $\sigma$ by
\begin{equation*}
\mu(x,1)~=~\left(\begin{array}{c} 0\\ -1 \end{array}\right) ,\quad  \mu(x,2)~=~\left(\begin{array}{c} 0\\ 0 \end{array}\right)
\quad\mbox{ and }\quad 
\sigma(x,1)~=~\sigma (x,2)= \left(\begin{array}{cc} 
0 &  0\\
0 &  0 
\end{array}\right)
\end{equation*}
for all $x\in\R^2$. Define the gain coefficients $g$ and $f$and the cost functions $c(.,1,2)$ and $c(.,2,1)$ by
\beqs
g(x,1)~=~g(x,2)~=~0\;, ~  f(x,1)~=~f(x,2) ~=~1 & \mbox{and} & c(x,1,2) ~ = ~ c(x,2,1)~=~c~>~0\;,
\enqs
for all $x\in\R^2$. Since the reward coefficients $f$ and $g$ do not depend on the state position $x$ we only need to focus on the constraint. In particular a strategy is optimal if it minimizes the number of switching orders and satisfy the state constraint.

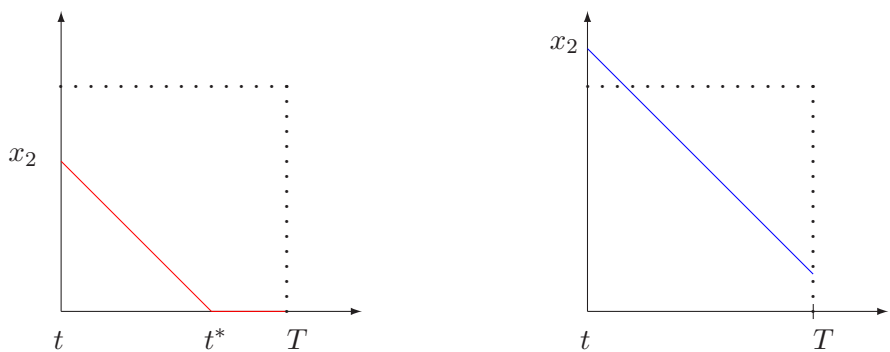
\begin{figure}[h]
\setlength{\unitlength}{1cm}

\begin{picture}(12,4.1)

\put(2,0){\vector(5,0){4}}
\put(2,0){\vector(0,5){4}}

\textcolor{red}{\put(2,2){\line(1,-1){2}}
\put(4,0){\line(1,0){1}}}
\put(1.3,2){$x_2$}

\put(5,0){\dottedline[.]{0.2}(0,0)(0,3)}
\put(5,0){\dottedline[.]{0.2}(-3,3)(0,3)}

\put(5,-0.5){$T$}
\put(1.9,-0.5){$t$}
\put(3.9,-0.5){$t^*$}

\put(9,0){\vector(5,0){4}}
\put(9,0){\vector(0,5){4}}
\put(12,-0.1){\line(0,1){0.2}}
\put(12,-0.5){$T$}
\put(8.9,-0.5){$t$}
\textcolor{blue}{\put(9,3.5){\line(1,-1){3}}}

\put(12,0){\dottedline[.]{0.2}(0,0)(0,3)}
\put(12,0){\dottedline[.]{0.2}(-3,3)(0,3)}
\put(8.5,3.5){$x_2$}

\end{picture}

\vspace{3mm}
\begin{footnotesize}\caption{\label{diagram-opt}Second component of optimal trajectories in the cases $x_2<T-t$,  $i=1$ (red curve) and $x_2\geq T-t$, $i=1$ (blue curve) .}
\end{footnotesize}

\end{figure}

As shown by Figure \ref{diagram-opt}, in the case $x_2<T-t$ and $i=1$,  the agent has to act at time $t^*$ to keep the second component non-negative (see the red curve). On the contrary, in the case $x_2\geq T-t$ and $i=1$, the blue curve shows that the system will satisfy the constraint until terminal time $T$ and there is no need to switch. We therefore get the following expression for the value function 
\begin{equation}\label{expr-v-ce}
v(t,x,1) ~ = ~ \left\{
\begin{array}{ccc}
T-t & \mbox{ if } & x_2~\geq~ T-t\;,\\
T-t-c & \mbox{ if } & x_2~ <~ T-t \;,
\end{array}
\right.
\end{equation}
for all $x=\left(\begin{array}{c}x_1 \\x_2\end{array}\right)\in\Dc$ and all $t\in[0,T]$. 

In particular the function $v(.,1)$ is discontinuous at
each point $(t,(x_1,T-t))$ for all $t\in[0,T]$ and all $x_1\in\R$. Hence the function $v$ is discontinuous even on the interior $Int(\Dc)$ of the constraint domain.
These discontinuities are induced by the state constraints that forces the operator to act so as to keep the diffusion in $\Dc$, even if this action is sub-optimal.

\section{Examples of application}\label{example elec}
We present in this section some models involving an optimal switching  problem under state constraint.
 \subsection{Hydroelectric pumped 
storage model}
The following simplified  hydroelectric pumped storage model is  inspired by \cite{carlud05}.    
Pumped Storage (currently, the dominant type of electricity storage) consists of large reservoir of water held by a hydroelectric dam at a higher elevation. When desired, the dam can be opened which activates the turbines and moves the water to another, lower reservoir. The generated electricity is sold to a power grid. As the water flows, the upper reservoir is depleted. Conversely, in times of low electricity demand, the water can be pumped back into the reservoir with required energy purchased from grid. 
A strategy $\alpha$ consists in a sequence of $\F$-stopping times $(\tau_{k})_{k}$ representing the intervention times and a sequence of $\Fc_{\tau_k}$-measurable random variables $(\zeta_{k})_{k}$ representing the changes of regime.
There are three possible regimes. 
\begin{enumerate}[(i)]
\item $\zeta_{k}$ $=$ $1$ : pump, in this case we set $\mu_1(x,1)=1$ and $\sigma_1(x,1)=0$.  
\item  $\zeta_{k}$ $=$ $2$: store, in this case we set $\mu_1(x,1)=0$ and $\sigma_1(x,1)=0$.
\item  $\zeta_{k}$ $=$ $3$: generate, in this case we set $\mu_1(x,1)=-1$ and $\sigma_1(x,1)=0$.
\end{enumerate}
For a given strategy $\alpha=(\tau_{k},\zeta_{k})_{k}$, we
denote by $L^\alpha_{t}$ the controlled water level in the upper reservoir. It satisfies the equation
\beqs\label{constI}
L^{\alpha}_{t} & = & L_0+\int_0^t \mu_1(L^\alpha_s,\alpha_s) ds+\int_0^t\sigma_1(L^\alpha_s,\alpha_s)dW_s\;, \qquad t\geq0\,.
\enqs
Denote by $P$ the electricity price process and suppose that it is a diffusion defined on $(\Omega,\Gc,\P)$ by
\beqs
P_t & = & P_0+\int_0^t\mu_2(P_s)ds+\int_0^t\sigma_2(P_s)dW_s\;,\quad t\geq 0\;.
\enqs
Let $X^\alpha$ be the controlled process defined by $X^\alpha=\begin{footnotesize}\left(\begin{array}{c}L^\alpha\\P\end{array}\right)\end{footnotesize}$. Then it satisfies the SDE
\beqs
X^\alpha_t & = & X_0+\int_0^t\mu(X_s^\alpha,\alpha_s)ds+\int_0^t\sigma(X_s^\alpha,\alpha_s)dW_s\;,\quad t\geq 0\;,
\enqs
with $\mu=\begin{footnotesize}\left(\begin{array}{c}\mu_1\\\mu_2\end{array}\right)\end{footnotesize}$ and $\sigma=\begin{footnotesize}\left(\begin{array}{c}\sigma_1\\\sigma_2\end{array}\right)\end{footnotesize}$.
 Suppose also that the cost of changing the regime from $i$ to $j$ is given by a constant $c(i,j)$. The expected pay-off for a strategy $\alpha$ is then given by
\beqs
J(0,X_0,\alpha) & = & \E\Big[\int_{0}^T-P_{t}dL^{\alpha}_{t}-\sum_{\tau_{k}\leq T}c(\zeta_{k-1},\zeta_{k})\Big]
 ~ = ~ \E\Big[\int_{0}^Tf(X^\alpha_t,\alpha_t)d{t}-\sum_{\tau_{k}\leq T}c(\zeta_{k-1},\zeta_{k})\Big]
\enqs
where $f$ is defined by $f(p,\ell,i)=-p\times \mu_1(\ell,i)$ for all $(p,\ell,i)\in\R\times\R\times\{1,2,3\}$.

Since the reservoir capacity is not infinite, the strategy $\alpha$ has to satisfy the constraint $0\leq L_{t}^\alpha\leq \ell_{max}$ for all $t\in[0,T]$. This corresponds to the general constraint $X^\alpha_t\in \Dc$ where $\Dc = \R\times[0,\ell_{max}]$. The goal of the energy producer is to maximize $J(0,X_0,\alpha)$ over the strategies $\alpha$ satisfying the constraint on the water level $L^\alpha$.

\subsection{Valuation of natural ressources}
The following model comes from \cite{BS85}. We consider an agent that holds a mine that produces a single homogeneous commodity. We suppose that the commodity price $S$ is given by
\beqs
S_t & = & S_0+\int_0^t\mu_1(S_u) du + \int_0^t\sigma_1(S_u) dW_u\;,\quad t\geq 0\;.
\enqs
The agent can choose to extract or not the commodity from the mine. Thus, the strategy $\alpha$ consists in a sequence of $\F$-stopping times $(\tau_{k})_{k}$ representing the intervention times and a sequence of $\Fc_{\tau_k}$-measurable random variables $(\zeta_{k})_{k}$ representing the changes of regime.
There are two possible regimes. 
\begin{enumerate}[(i)]
\item $\zeta_{k}$ $=$ $1$: extraction, in this case we set $\mu_2(x,1)=-1$ and $\sigma_2(x,1)=0$.  
\item  $\zeta_{k}$ $=$ $0$: no extraction, in this case we set $\mu_2(x,2)=0$ and $\sigma_2(x,2)=0$.
\end{enumerate}
For a strategy $\alpha=(\tau_{k},\zeta_{k})_{k}$, we denote by $Q^\alpha_t$  the physical inventory of the mine at time $t$. Therefore, it satisfies the equation
\beqs
Q^{\alpha}_{t} & = & Q_0+\int_0^t \mu_2(Q^\alpha_s,\alpha_s) ds+\int_0^t\sigma_2(Q^\alpha_s,\alpha_s)dW_s\;, \qquad t\geq0\,.
\enqs
 Denote by $X^\alpha$ the controlled process defined by $X^\alpha=\begin{footnotesize}\left(\begin{array}{c}S\\Q^\alpha\end{array}\right)\end{footnotesize}$. Then it satisfies the SDE
\beqs
X^\alpha_t & = & X_0+\int_0^t\mu(X_s^\alpha,\alpha_s)ds+\int_0^t\sigma(X_s^\alpha,\alpha_s)dW_s\;, \qquad t\geq0\,,
\enqs
with $\mu=\begin{footnotesize}\left(\begin{array}{c}\mu_1\\\mu_2\end{array}\right)\end{footnotesize}$ and $\sigma=\begin{footnotesize}\left(\begin{array}{c}\sigma_1\\\sigma_2\end{array}\right)\end{footnotesize}$.
 Suppose also that the cost of changing the regime from $i$ to $j$ is given by a constant $c(i,j)$. The expected pay-off for a strategy $\alpha$ is then given by
\beqs
J(0,X_0,\alpha) & = & \E\Big[\int_{0}^TS_{t}dQ^{\alpha}_{t}-\sum_{\tau_{k}\leq T}c(\zeta_{k-1},\zeta_{k})\Big]
 ~~ = ~~ \E\Big[\int_{0}^Tf(X^\alpha_t,\alpha_t)d{t}-\sum_{\tau_{k}\leq T}c(\zeta_{k-1},\zeta_{k})\Big]
\enqs
where $f$ is defined by $f(s,q,i)=-s\times \mu_2(q,i)$ for all $(s,q,i)\in\R\times\R\times\{0,1\}$.

Since the physical inventory is non-negative, the strategy $\alpha$ has to satisfy the constraint $ Q_{t}^\alpha\geq 0$ for all $t\in[0,T]$. This corresponds to the general constraint $X^\alpha_t\in \Dc$ where $\Dc = \R\times\R_+$. Thus, the aim of the agent is to maximize $J(0,X_0,\alpha)$ over the strategies $\alpha$ satisfying the constraint on the inventory $Q^\alpha$.

\subsection{Reversible technology investment}

We present a simplified version of the model studied in \cite{lpv}. We consider a firm whose activities generate cash process by using some technology.  The firm has at any time the possibility to choose between two
technologies: a modern one and an old one. Therefore,  its strategy $\alpha$ consists in a sequence of $\F$-stopping times $(\tau_{k})_{k}$ representing the times of change of technology and a sequence of $\Fc_{\tau_k}$-measurable random variables $(\zeta_{k})_{k}$ representing the chosen technology at each time $\tau_k$.
Thus, there are two possible regimes. 
\begin{enumerate}[(i)]
\item $\zeta_{k}$ $=$ $1$: old technology, in this case we set $\mu(x,1)=\delta_1 x $ and $\sigma(x,1)=\gamma_1 x$.  
\item  $\zeta_{k}$ $=$ $2$: modern technology, in this case we set $\mu(x,2)=\delta_2 x $ and $\sigma(x,2)=\gamma_2 x$.
\end{enumerate}
Here $\gamma_1$, $\gamma_2$, $\delta_1$ and $\delta_2$ are four constants with $\delta_1<\delta_2$ and $\gamma_1<\gamma_2$ (the modern technology has a better rate but a worse uncertainty than the old technology). For a strategy $\alpha=(\tau_{k},\zeta_{k})_{k}$, we denote by $X^\alpha_t$ the cash reserve at time $t$ of the firm. We suppose that it satisfies the equation
\beqs
X^{\alpha}_{t} & = & X_0+\int_0^t \mu(X^\alpha_s,\alpha_s) ds+\int_0^t\sigma(X^\alpha_s,\alpha_s)dW_s\;, \qquad t\geq0\,.
\enqs
We also suppose that the cost of changing the technology from $i$ to $j$ is given by a constant $c(i,j)$. Then the expected pay-off at terminal time $T$ for a strategy $\alpha$ is given by
\beqs
J(0,X_0,\alpha) & = & \E\Big[X^{\alpha}_{T}-\sum_{\tau_{k}\leq T}c(\zeta_{k-1},\zeta_{k})\Big]
\enqs
We suppose that the firm have to satisfy the following solvency constraint 
$ X_{t}^\alpha\geq 0$ for all $t\in[0,T]$. This corresponds to the constraint domain $\Dc = \R_+$. Thus, the goal of the firm is to maximize $J(0,X_0,\alpha)$ over the strategies $\alpha$ satisfying the constraint on the cash reserve $R^\alpha$.

\section{Unconstrained penalized switching problem}
\setcounter{equation}{0} \setcounter{Assumption}{0}
\setcounter{Theorem}{0} \setcounter{Proposition}{0}
\setcounter{Corollary}{0} \setcounter{Lemma}{0}
\setcounter{Definition}{0} \setcounter{Remark}{0}
\subsection{An unconstrained penalized approximating problem}

We now introduce an 
approximation of our initial constrained problem.
 This approximation consists in a penalization of the coefficients $f$ and $g$ out of the domain $\Dc$ where the controlled underlying diffusion is constrained to stay.
 
Consider, for $n\geq1$, 
the functions
$f_{n}:\R^d\times\Ic\rightarrow\R$ and $g_{n}:\R^d\times\Ic\rightarrow\R$ 
defined by
\beq
f_{n}(x,i) & = & 
f(x,i)-n\Theta_{n}(x)\;,\\
g_{n}(x,i) & = & 
g(x,i)-n\Theta_{n}(x)\;,
\enq
for all $(x,i)\in\R^d\times\Ic$, where the function $\Theta_{n}:\R^d\rightarrow[0,1]$ is given by
\beq
\Theta_{n}(x) & = & n\Big(d\big(x,\Dc\big)\wedge\frac{1}{n}\Big)
~=~ nd(x,\Dc)\wedge1\,,
\enq
with
 $d(x,\Dc)$ $=$ $\inf_{x'\in\Dc}|x-x'|$ for all $x\in\R^d$.

 Given an initial condition $(t,x)$ and a switching control $\alpha=(\tau_k,\zeta_k)_{k\geq0}\in\Ac_t$, we consider the
total penalized profit starting from
$(t,x,i)\in[0,T]\times\R^{d}\times\Ic$ at horizon $T$, defined
by:
\beqs
J_{n}(t,x,\alpha) & = & \E\Big[g_n\big(X^{t,x,\alpha}_T,\alpha_T\big)+\int_t^Tf_n\big(X^{t,x,\alpha}_s,\alpha_s\big)ds-\sum_{k\geq 1}c\big(X^{t,x,\alpha}_{\tau_k},\zeta_{k-1},\zeta_{k}\big)\mathds{1}_{\tau_k\leq T}\Big]\;.
\enqs
We can then define the penalized unconstrained  value function $v_{n}:~[0,T]\times\R^{d}\times\Ic\rightarrow\R$ by
\beq
v_{n}(t,x,i) & = & \sup_{\alpha\in\Ac_{t,i}} J_{n}(t,x,\alpha)\;,
\enq
for all $n\geq 1$ and all $(t,x,i)\in[0,T]\times\R^d\times\Ic$, where $\Ac_{t,i}$ is the set of strategies $\alpha=(\tau_k,\zeta_k)_{k\geq0}\in\Ac_t$ such that $\zeta_0=i$.
\subsection{Convergence of the penalized unconstrained problems}
We now state the main result of this section which concerns the convergence of the functions $v_n$ to $v$. The main line of the proof is to take a sequence of almost optimal strategies for the functions $v_n$ and to make it converge to a strategy that we expect to be optimal. To do this, we need to prove measurability and convergence results for diffusion driven by a converging sequence of Brownian motions. These results are presented in details in the Appendix \ref{apptd. means cvgce}.
\begin{Theorem}\label{Propconv}
Under \textbf{(H1)} and \textbf{(H2)}, the sequence $(v_{n})_{n\geq1}$ is nonincreasing and converges
 on $[0,T]\times\Dc \times\Ic$ to the function $v$:
\beq
v_{n}(t,x,i) & \downarrow & v(t,x,i) \qquad \mbox{as }~~ n\uparrow +\infty, 
\enq 
for all $(t,x,i)\in [0,T]\times\Dc\times\Ic$. Moreover, for any 
 $(t,x,i)\in [0,T]\times\Dc\times\Ic$,
there exists a strategy  $\alpha^*\in\Ac_{t,x,i}^\Dc$ such that
\beqs
v(t,x,i) & = &  
J(t,x,\alpha^*)\;.
\enqs 
\end{Theorem}

\ni \textbf{Proof.}
Fix $(t,x,i)\in[0,T]\times\Dc\times\Ic$. 
Since $f_{n+1}\leq f_n$ and $g_{n+1}\leq g_n$ we get 
\beqs
J_{n+1}(t,x,\alpha) & \leq & J_n(t,x,\alpha)\;, 
\enqs
 for all $n\geq1$ and $\alpha\in\Ac_t$. From this last inequality we deduce that
\beqs
v_{n+1}(t,x,i) & \leq & v_{n}(t,x,i)\;,\quad n\geq1\;.
\enqs
 We now prove that $(v_n)_n$ converges to $v$.  We first notice that 
 \beqs
 J_n(t,x,\alpha) & = & J(t,x,\alpha)\;,
 \enqs
 for any $n\geq1$, any initial condition $(t,x,i)\in[0,T]\times\Dc\times\Ic$ and any switching strategy $\alpha\in \Ac_{t,x,i}^{\Dc}$. Therefore, we get $v_n\geq v$ for all $n\geq 1$. Denote by $\bar v$ the pointwise limit of $(v_n)_n$:
\beqs
\bar v (t,x,i) & = & \lim_{n\rightarrow\infty}v_{n}(t,x,i)\;, \quad (t,x,i)\in[0,T]\times\Dc\times\Ic\;.
\enqs
 Then we have $\bar v(t,x,i)\geq v(t,x,i)$. If $\bar v(t,x,i)=-\infty$  we obviously get $\bar v(t,x,i)= v(t,x,i)$.
 
 We now suppose that $\bar v(t,x,i)>-\infty$ and prove that $\bar v(t,x,i)\leq v(t,x,i)$. 
 We proceed in $3$ steps.
 \vspace{2mm}
 

 \noindent\textbf{Step 1.} \textit{Convergence of a sequence of almost optimal strategies for the unconstrained problems.}

 \ni \textbf{Substep 1.1.} \textit{Bounded sequence of almost optimal strategies.} 
 
 \ni For $n\geq 1$, let  $\alpha^n=(\tau_k^n,\zeta_k^n)_{k\geq0}\in \Ac_{t,i}$ a switching strategy such that
 \beqs
 J_n(t,x,\alpha^n) & \geq & v_n(t,x,i)-{1\over n}\;.
 \enqs
 We can suppose without loss of generality that 
 \beq\label{bornetaunk}
 \tau_k^n\in [0,T]\cup\{T+1\}\qquad\P-a.s.
 \enq
 for all $n\geq 1$ and all $k\geq 0$. Indeed, fix $n\geq 1$ and consider the strategy $\hat \alpha^n=(\hat \tau_k^n,\hat \zeta_k^n)_{k\geq0}\in \Ac_{t,i}$ defined by 
\beqs
\hat \tau_k^n & = & \tau_k^n\mathds{1}_{\tau_k^n\leq T} +(T+1)\mathds{1}_{\tau_k^n>T}\;, \\
\hat \zeta_k^n & = & \zeta_k^n\mathds{1}_{\tau_k^n\leq T}+ i\mathds{1}_{\tau_k^n> T}\;.
\enqs 
 Then we have  $J_n(t,x,\alpha^n)=J_n(t,x,\hat \alpha^n)$ and we can replace $\alpha^n$ by $\hat \alpha^n$ which satisfies \reff{bornetaunk}.
 
 \vspace{2mm}
 
\ni  \textbf{Substep 1.2.} \textit{Tightness  and convergence of $(W,\alpha^n)_n$.}

\ni We now prove that the sequence of $C([0,T],\R^d)\times \big(\R_+\times \Ic\big)^\N$-valued random variables $(W,\alpha^n)_{n\geq 1}$ is tight.
 Fix a sequence $(\delta_\ell)_\ell$ of positive numbers such that
 \beq\label{cond delta}
 \delta_\ell~\xrightarrow[\ell\rightarrow\infty]{}  ~0 & \mbox{ and } & 2^\ell\delta_\ell\ln\big({2T \over \delta_\ell}\big) ~\xrightarrow[\ell\rightarrow\infty]{}  ~0\;.
 \enq
We define for $\eta>0$ and $C>0$ the subset $\Kc_\eta^C$ of $C([0,T],\R^d)$ by
 \beqs
 \Kc_\eta^C & = &\bigcap_{\ell\geq1} \Kc_{\eta,\ell}^C 
 \enqs
 where
 \beqs
 \Kc_{\eta,\ell}^C & = & \left\{h\in C([0,T],\R^d)~:~h(0)=0 ~\mbox{ and }~{\rm mc}_{\delta_\ell}(h)~\leq~C{2^\ell\delta_\ell\ln\big({2T \over \delta_\ell}\big)\over \eta}
 \right\}
 \enqs
 and ${\rm mc}$ denotes the modulus of continuity defined by
 \beqs
{\rm mc}_\delta(h) & = & \sup_{\begin{tiny}\begin{array}{c}s,t\in[0,T]\\|s-t|\leq \delta\end{array}
 \end{tiny}}\big|h(s)-h(t)\big|
 \enqs
 for any $h\in C([0,T],\R^d)$ and any $\delta>0$. Using Arz\'ela-Ascoli theorem, we get from \reff{cond delta} that $\Kc^C_\eta$ is a compact subset of $C([0,T],\R^d)$. 
 We now define the subset $\mathbf{K}^C_\eta$ of $C([0,T],\R^d)\times\big(\R_+\times\Ic\big)^\N$ by
 \beqs
\mathbf{K}^C_\eta  & = & \Kc^C_\eta\times \big([0,T+1]\times\Ic\big)^\N\;.
 \enqs
 From Tychonov theorem and since $\Kc^C_\eta$ is compact, we get that  $\mathbf{K}^C_\eta$ is a compact subset of $C([0,T],\R^d)\times\big(\R_+\times\Ic\big)^\N$ endowed with the norm $\|\cdot\|$ defined by 
 \beqs
 \big\|\big(h,{(t_k,z_k)}_{k\geq0}\big)\big\| & = & \sup_{t\in[0,T]}|h(t)|+ \sum_{k\geq0}{(|t_k|+|z_k|)\wedge 1\over 2^k}
 \enqs
  for all $h\in C([0,T],\R^d)$ and ${(t_k,z_k)}_{k\geq0}\in (\R_+\times\Ic)^\N$.
 We then have from \reff{bornetaunk} 
 \beqs
 \P\Big((W,\alpha^n)\in \mathbf{K}^C_\eta\Big) & = &  \P\Big(W\in \Kc^C_\eta\Big)
 \enqs  
 for all $\eta>0$, $C>0$ and $n\geq1$. Using Markov inequality we get
 \beq\nonumber
 \P\Big(W\in \Kc^C_\eta\Big) & = & 1-\P\Big(W\notin \Kc^C_\eta\Big)\\ \nonumber
  & \geq & 1-\sum_{\ell\geq 1}\P\Big(W\notin \Kc^C_{\eta,\ell}\Big)\\ \label{ineqMarkov}
  & \geq & 1-\sum_{\ell\geq 1}{\E\Big[{\rm mc}_{\delta_\ell}(W)\Big]\over C{2^\ell\delta_\ell\ln\big({2T \over \delta_\ell}\big)\over \eta}}\;.
 \enq
 From Theorem 1 in \cite{FN08}, there exists a constant $C^*$ such that
 \beq\label{ineqFN}
 \E\Big[{\rm mc}_\delta(W)\Big] & \leq & C^*\delta\ln\big({2T \over \delta}\big)\;.
 \enq
 for all $\delta>0$.  Therefore, we get from \reff{ineqMarkov} and \reff{ineqFN}
  \beqs
  \P\Big((W,\alpha^n)\in \mathbf{K}^{C^*}_\eta\Big) & \geq & 1-\eta \;,
  \enqs
for all $\eta\in(0,1)$, and the sequence $(W,\alpha^n)_n$ is tight. 

We deduce from Prokhorov theorem that, up to a subsequence, 
\beqs
\P\circ(W,\alpha^n)^{-1} & \xrightarrow[n\rightarrow\infty]{} &\Lc.
\enqs
with $\Lc$ a probability measure on $\big(C([0,T],\R^d)\times(\R\times\Ic)^\N,\|\cdot\|\big)$.
\vspace{2mm}

 \noindent\textbf{Step 2.} \textit{Change of probability space.}

\ni 
 Since $\big(C([0,T],\R^d)\times(\R\times\Ic)^\N,\|\cdot\|\big)$ is separable, we get from the Skorokhod representation theorem 
that there exists a probability space $(\tilde \Omega,\tilde \Gc,\tilde \P)$ on which are defined Brownian motions $\tilde W^{n}$, $n\geq 1$, and $\tilde W$, and  random variables $\tilde \alpha^{n}=(\tilde \tau_k^{n},\tilde \zeta_k^{n})_{k\geq 0}$, $n\geq 1$, and $\tilde \alpha=(\tilde \tau_k^{},\tilde \zeta_k^{})_{k\geq 0}$ 
 such that
\beq\label{egal-loi-n-ntilde}
\tilde \P\circ ( \tilde W^{n}, \tilde \alpha^{n})^{-1} & = & \P \circ ( W, \alpha^n)^{-1}
\enq
for all $n\geq 1$ and 
\beq\label{convpsW*alpha*}
\Big\|\big( \tilde W^n, \tilde \alpha^{n} \big)-\big(\tilde W, \tilde  \alpha\big)\Big\| &\xrightarrow[n\rightarrow\infty]{\tilde \P-a.s.}  & 0\;.
\enq 
In particular we get 
\beqs
\Lc & = & \tilde \P \circ (\tilde W, \tilde \alpha)^{-1}\;.
\enqs

\vspace{2mm}

 \noindent\textbf{Substep 2.1} \textit{Measurability properties for $\tilde\alpha^n$ and $\tilde\alpha$.} 
 


\ni We now prove that each $\tilde \tau_k$ is an $\tilde \F$-stopping time and $\zeta_k$ is $\tilde \Fc_{\tilde \tau_k}$-measurable where $\tilde \F=(\tilde \Fc_t)_{t\geq 0}$ is the complete right-continuous filtration generated by $\tilde W$.

For $n\geq 1$, denote by $\tilde \F^n=(\tilde \Fc^n_t)_{t\geq0}$ the complete right-continuous filtration generated by $\tilde W^n$. Using Proposition \ref{egal-loi-meas}, we get from \reff{egal-loi-n-ntilde} that $\tilde \tau_k^n$ is an $\tilde \F^n$-stopping time and that $\tilde \zeta^n_k$ is $\tilde \Fc^n_{\tilde \tau^n_k}$-measurable for all $n\geq 1$ and $k\geq 0$. Then using Proposition \ref{prop-meas-conv}, we get from \reff{convpsW*alpha*} that $\tilde \tau_k$ is an $\tilde \F$-stopping time and that $\tilde \zeta_k$ is $\tilde \Fc_{\tilde \tau_k}$-measurable for all $k\geq 0$.

\vspace{2mm}

\vspace{2mm}

\ni \textbf{Substep 2.2.} \textit{Equality of the penalized gains  and convergence of the associated controlled diffusions.}

\ni From the previous substep, we can define the diffusions $\tilde X ^{t,x,\tilde \alpha^n}$ and $\tilde X^{t,x,\tilde \alpha}$ on $(\tilde \Omega,\tilde \Gc,\tilde \P)$ by
\beqs
\tilde X ^{t,x,\tilde \alpha^n}_s & =& x+\int_t^sb(\tilde X^{t,x,\tilde \alpha^n} _r,\tilde  \alpha ^n _r)dr+\int_t^s\sigma(\tilde X^{t,x,\tilde \alpha^n}_r,\tilde  \alpha ^n _r)d\tilde W ^n_r\;,\quad s\geq t,
\enqs
and 
\beqs
\tilde X^{t,x,\tilde \alpha} _s & =& x+\int_t^sb(\tilde X^{t,x,\tilde \alpha}_r,\tilde  \alpha  _r)dr+\int_t^s\sigma(\tilde X_r^{t,x,\tilde \alpha},\tilde  \alpha  _r)d\tilde W_r\;,\quad s\geq t,
\enqs
and the associated gains $J_n(t,x,\tilde \alpha^n)$ and $J(t,x,\tilde \alpha)$ by
\beqs
\tilde J_n (t,x,\tilde \alpha ^n) & = & \E^{\tilde \P}\Big[g_n\big(\tilde X^{t,x,\tilde \alpha^n}_T,\tilde \alpha^n_T\big)+\int_t^Tf_n\big(\tilde X^{t,x,\tilde \alpha^n}_s,\tilde \alpha^n_s\big)ds-\sum_{k\geq 1}c\big(\tilde X^{t,x,\tilde \alpha^n}_{\tilde\tau^n_k},\tilde \zeta^n_{k-1},\tilde \zeta^n_{k}\big)\mathds{1}_{\tilde\tau^n_k< T}\Big]
\enqs
and 
\beqs
\tilde J (t,x,\tilde \alpha ) & = & \E^{\tilde \P}\Big[g\big(\tilde X^{t,x,\tilde \alpha}_T,\tilde \alpha_T\big)+\int_t^Tf\big(\tilde X^{t,x,\tilde \alpha}_s,\tilde \alpha_s\big)ds-\sum_{k\geq 1}c\big(\tilde X^{t,x,\tilde \alpha}_{\tilde\tau_k},\tilde \zeta_{k-1},\tilde \zeta_{k}\big)\mathds{1}_{\tilde\tau_k< T}\Big]\;.
\enqs
Since $(W,\alpha^{n})$ and $(\tilde W^n, \tilde \alpha^n)$ have the same law, we deduce from \textbf{(H1)} and \textbf{(H2)} that
\beq\label{eq-pen-gain}
J_n(t,x,\alpha^{n}) & = & \tilde J_n(t,x, \tilde  \alpha^n)~~\geq~~v_n(t,x,i)-{1\over n}\;,\quad n\geq 1\;.
\enq

\vspace{2mm}

%
%

\vspace{2mm}

  
\ni We now prove that, up to a subsequence, 
\beq\label{majlimsup}
\limsup_{n\rightarrow\infty}\tilde J_n( t,x, \tilde\alpha^{n}) & \leq & \tilde J(t,x,\tilde \alpha).
\enq 
We first notice that
$\limsup_{n\rightarrow\infty}\tilde J_n( t,x, \tilde \alpha^{n})  \leq   \limsup_{n\rightarrow\infty}\tilde J( t,x,\tilde \alpha^{n})$. 
From Proposition \ref{AnConvDiff} and \reff{convpsW*alpha*} we have
\beq\label{convDC}
\E^{\tilde \P} \Big[ \sup_{s\in[t,T]}\big|\tilde X^{t,x,\tilde \alpha}_s-\tilde X_s^{t,x,\tilde \alpha^n} \big|^2 \Big] & \xrightarrow[n\rightarrow\infty]{} & 0\;.
\enq
We therefore get, up to a subsequence,
\beq\label{convpsXnXinfty}
\sup_{s\in[t,T]}\left|\tilde X^{t,x,\tilde \alpha^{n}}_{s}-\tilde X^{t,x,\tilde \alpha^{}}_{s}\right|  & \xrightarrow[n\rightarrow\infty]{\tilde \P-a.s.}  & 0\;.
\enq
This implies with \textbf{(H2)} (i) and (ii) and \reff{convpsW*alpha*}
\beqs
g\big(\tilde X^{t,x,\tilde \alpha^{n}}_T,\tilde \alpha^{n}_T\big)+\int_t^Tf\big(\tilde X^{t,x,\tilde \alpha^{n}}_s,\tilde \alpha^{n}_s\big)ds & \xrightarrow[n\rightarrow\infty]{\tilde \P-a.s.}  & 
g\big(\tilde X^{t,x,\tilde \alpha}_T,\tilde \alpha_T\big)+\int_t^Tf\big(\tilde X^{t,x,\tilde \alpha}_s,\tilde \alpha_s\big)ds\;. 
\enqs
Moreover, since $\bar v (t,x,i)>-\infty$ we have from \textbf{(H2)} (ii)
\beqs
\sup_{n\geq 1}\#\big\{k\geq1~:~\tilde \tau_k^{n}\leq T\big\} & < & +\infty\;,\quad\tilde \P-a.s.
\enqs
This last estimate, \reff{bornetaunk}, \reff{convpsW*alpha*} and \reff{convpsXnXinfty} imply
\beqs
\liminf_{n\rightarrow\infty}\sum_{k\geq 1}c\big(\tilde X^{t,x,\tilde \alpha^{n}}_{\tilde \tau_k},\tilde \zeta^{n}_{k-1},\tilde \zeta^{n}_{k}\big)\mathds{1}_{\tilde \tau_k^{n}\leq T} & 
\geq &  \sum_{k\geq 1}c\big(\tilde X^{,t,x,\tilde \alpha}_{\tilde \tau_k},\tilde \zeta_{k-1},\tilde \zeta_{k}\big)\mathds{1}_{\tilde \tau_k\leq T}\;,\quad \tilde \P-a.s.
\enqs
We finaly conclude by using Fatou's Lemma. 

 \vspace{2mm}

\ni \textbf{Substep 2.3}  \textit{The process $\tilde X^{t,x,\tilde \alpha}$ satisfies the constraint $\tilde X^{t,x,\tilde \alpha} _s\in \Dc$ for all $s\in [t,T]$.}

\ni 
For $\eps>0$, we define the set $\Dc_\eps$ by
\beqs
\Dc_\eps & = & \Big\{ x'\in\R^d~:~d(x',\Dc)< \eps \Big\}\;.
\enqs
Suppose that there exists some $\eps>0$ such that 
\beqs
\E^{\tilde \P}\Big[ \int_t^T\mathds{1}_{\Dc_\eps^c}(\tilde X_{s}^{t,x,\tilde \alpha})ds \Big] & > & 0\;.
\enqs
From \reff{convpsXnXinfty} and the dominated convergence theorem we can find $\eta>0$ and $n_\eta\geq 1$ such that, up to a subsequence,
\beqs
\E^{\tilde \P}\Big[ \int_t^T\mathds{1}_{\Dc_\eps^c}(\tilde X_{s}^{t,x,\tilde \alpha^{n}})ds \Big] & \geq & \eta
\enqs
for all $n\geq n_\eta$.
From the definition of $f_n$ and $g_n$ and the previous inequality, there exists a constant $C$ such that
\beqs
\tilde J(t,x,\tilde \alpha^{n}) & \leq & C\E^{\tilde \P}\Big[\sup_{s\in[t,T]}\big|\tilde X^{t,x,\tilde \alpha^{n}}_{s}\big|\Big]-n\eta 
\enqs
for any $n\geq {1\over \eps}\vee n_\eta$. Sending $n$ to infinity we get from \reff{eq-pen-gain} and \reff{XS2} applied on $(\tilde \Omega,\tilde \Gc,\tilde \P)$ 
\beqs
\bar v (t,x,i) & = & \lim_{n\rightarrow\infty}\tilde J_n(t,x,\tilde \alpha^{n}) ~~=~~-\infty
\enqs
which contradicts $\bar v (t,x,i)>-\infty$.
We therefore obtain
\beqs
\E^{\tilde \P}\Big[\int_t^T\mathds{1}_{\Dc_\eps^c}(\tilde X^{t,x,\tilde \alpha}_s)ds\Big] & = & 0
\enqs
for all $\eps>0$ and $\E^{\tilde \P}\Big[\int_t^T\mathds{1}_{\{\tilde X^{t,x,\tilde \alpha}_s\notin\Dc\}}ds\Big]  =  0$.
 Since $\tilde X^{t,x,\tilde \alpha}$ is continuous, we get 
 \beqs
 \tilde \P\Big(\tilde X^{t,x,\tilde \alpha}\in \Dc\;,~\forall s\in[t,T]\Big) & = & 1\;.
 \enqs

\vspace{2mm}

 \noindent\textbf{Step 3.} \textit{Back to $(\Omega,\Gc,\P)$ and conclusion.}
 
 
\ni We construct $\alpha^*\in \Ac_{t,i}$ such that $(W,\alpha^*)$ has the same law as $(\tilde W,\tilde \alpha)$. Using Proposition \ref{ecriture-ta} we can find Borel functions $\psi_k$ and $\phi_k$, $k\geq 1$ such that 
 \beqs
 \tilde \tau_k~~=~~\psi_k\big((\tilde W_s)_{s\in[0,T]}\big) & \mbox{ and } & \tilde \zeta_k~~=~~\phi_k\big((\tilde W_s)_{s\in[0,T+1]}\big)\quad \tilde \P-a.s.
 \enqs 
 for all $k\geq 0$. Define the strategy $\alpha^*=(\tau^*_k,\zeta^*_k)_{k\geq 0}$ by
 \beqs
 \tau^*_k~~=~~\psi_k\big(( W_s)_{s\in[0,T]}\big) & \mbox{ and } &  \zeta^*_k~~=~~\phi_k\big(( W_s)_{s\in[0,T+1]}\big)
 \enqs
 for all $k\geq 0$. Obviously $(W,\alpha^*)$ has the same law as $(\tilde W,\tilde \alpha)$. Moreover, from Proposition \ref{egal-loi-meas}, each $\tau_k^*$ is an $\F$-stopping time and each $\zeta_k^*$ is $\Fc_{\tau^*_k}$-measurable. We deduce that $\alpha^*\in \Ac_{t,i}$.
 Using Substep 2.3 we also get $\alpha^*\in \Ac_{t,x,i}^{\Dc}$.
\ni  From \reff{eq-pen-gain} and \reff{majlimsup} we get, up to a subsequence,
\beqs
\tilde J (t,x,\tilde \alpha)~~\geq~~\limsup_{n\rightarrow\infty}\tilde J_n(t,x,\tilde \alpha^{n})~~=~~\limsup_{n\rightarrow\infty} J_n(t,x,\alpha^{n}) & \geq & \bar v(t,x,i)\;.
\enqs
Since $(W,\alpha^*)$ and $(\tilde W,\tilde \alpha)$ have the same law and $\alpha^*\in\Ac_{t,x,i}^\Dc$ we get 
\beqs
v(t,x,i)~~\geq~~ J(t,x,\alpha^*) & = & \tilde J(t,x,\tilde \alpha)~~\geq~~ \bar v(t,x,i)\;.
\enqs
\ep

In general,  proving a regularity result on  the value function of a constrained optimization problem is very technical (see e.g. \cite{Son86} or \cite{Kat}). In our case, Theorem \ref{Propconv} gives  a semi-regularity for $v$. 
\begin{Corollary}\label{corvusc}
Under \textbf{(H1)} and  \textbf{(H2)}, the function $v(.,i)$ is upper semicontinuous on $[0,T)\times\Dc$ for all $i\in\Ic$.
\end{Corollary}
\ni\textbf{Proof.}
 Fix $i\in \Ic$. From \textbf{(H1)} and  \textbf{(H2)} the value function $v_n(.,i)$ associated to the penalized optimal switching problem is  continuous  on $[0,T)\times\R^d$ (see e.g. 
  \cite{bou06}). 
%
%
From Theorem \ref{Propconv}, the function $v(.,i)$ is upper semicontinuous on $[0,T)\times\Dc$ as an infimum of continuous functions. \ep

\section{Dynamic programming and variational inequalities}
\setcounter{equation}{0} \setcounter{Assumption}{0}
\setcounter{Theorem}{0} \setcounter{Proposition}{0}
\setcounter{Corollary}{0} \setcounter{Lemma}{0}
\setcounter{Definition}{0} \setcounter{Remark}{0}

\subsection{The dynamic programming principle}
In this section we state the dynamic programming principle.
We first need the following lemmata. We postpone their proofs to the Appendix \ref{appendice preuve lemme} to focus on the dynamic programming principle and its proof.
\begin{Lemma}\label{Propfngn}
Under \textbf{(H2)},  the functions $f_{n}$ and $g_{n}$ are locally Lipschitz continuous and have polynomial growth:
\begin{itemize}
\item  for any $n\geq 1$ and any $R>0$, there exists a constant $L_{R,n}$ such that
\beqs
|g_{n}(x,i)-g_{n}(x',i)|+|f_{n}(x,i)-f_{n}(x',i)| & \leq & L_{R,n}|x-x'|,
\enqs
for all $x,x'\in\R^{d}$ such that $|x|\leq R$ and  $|x'|\leq R$, and all $i\in\Ic$.
\item  for any $n\geq 1$, there exists a constant $C_n$ such that 
\beqs
|g_{n}(x,i)|+|f_{n}(x,i)| & \leq & C_n\big(1+|x|^q\big),
\enqs
for all $x\in\R^{d}$ and all $i\in\Ic$.
\end{itemize}
\end{Lemma}
\begin{Lemma}\label{prop-bound-v} Under \textbf{(H1)} and \textbf{(H2)}, there exists a constant $C$ such that 
\beq\label{maj-lin-v}
v_n(t,x,i) & \leq & C\big(1+|x|^q\big)
\enq
for all $n\geq 1$ and all $(t,x,i)\in[0,T]\times\Dc\times\Ic$.
\end{Lemma}

\vspace{2mm}

\ni We are now able to state the dynamic programming principle.
\begin{Theorem}\label{PPD}Under \textbf{(H1)} and \textbf{(H2)}, the value function $v$ satisfies the following dynamic
programming equality: 
\beq\label{DPP1}
v(t,x,i) & = & \sup_{\alpha=(\tau_{k},\zeta_{k})_{k}\in\Ac^{\Dc}_{t,x,i}}\E\Big[\int_{t}^{\nu}f(X^{t,x,\alpha}_{s},\alpha_{s})ds\nonumber-\sum_{t\leq\tau_{k}\leq\nu}c(X^{t,x,\alpha}_{\tau_k},\zeta_{k-1},\zeta_{k})\\
 & & \qquad\qquad\qquad\qquad\qquad\qquad\qquad\qquad\qquad+v\Big(\nu,X^{t,x,\alpha}_\nu,\alpha_{\nu}\Big)\Big]\;.\qquad
\enq
for any $(t,x,i)\in [0,T]\times\Dc \times \Ic$, and any stopping time $\nu$ valued in $[t,T]$. 
\end{Theorem}

\ni\textbf{Proof.} We first notice that the l.h.s. of \reff{DPP1} is well defined. Indeed, for a given stopping time $\nu$ valued in $[t,T]$ and a strategy $\alpha\in\Ac^{\Dc}_{t,x,i}$, we get from the regularity of $v$ given by Corollary \ref{corvusc} that the random quantity $v\Big(\nu,X_\nu^{t,x,\alpha},\alpha_{\nu}\Big)$ is measurable. Moreover, from Lemma  
\ref{prop-bound-v}, \reff{XS2} and the inequality $v\leq v_n$, we get that its expectation is well defined. 

Fix $(t,x,i)\in[0,T]\times\Dc\times\Ic$. If $\Ac^{\Dc}_{t,x,i}=\emptyset$ then the two hand sides of \reff{DPP1} are equal to $-\infty$ so the equality holds. 

Suppose now that $\Ac^{\Dc}_{t,x,i}\neq\emptyset$ and let $\alpha=(\tau_{k},\zeta_{k})_{k}\in\Ac^{\Dc}_{t,x,i}$ and $\nu$ a stopping time valued in $[t,T]$. 
From Lipschitz properties of $f_{n}$ and $g_{n}$ given by Lemma \ref{Propfngn}, we have by Lemma 4.4 in \cite{bou06}
\beq
v_{n}(t,x,i) & \geq & \E\Big[\int_{t}^{\nu}f_{n}\big(X^{t,x,\alpha}_{s},\alpha_{s}\big)ds\nonumber-\sum_{t\leq\tau_{k}\leq\nu}c\big(X^{t,x,\alpha}_{\tau_k},\zeta_{k-1},\zeta_{k}\big)
 +v_{n}\big(\nu,X_{\nu}^{t,x,\alpha},\alpha_{\nu}\big)\Big]\;,
\enq
for all $n\geq 1$. Since $\alpha\in\Ac^{\Dc}_{t,x,i}$ we have from the definition of $f_{n}$,   
\beqs
f_{n}(X^{t,x,\alpha}_{s},\alpha_{s}) & = & f(X^{t,x,\alpha}_{s},\alpha_{s})
\enqs 
for $d\P\otimes ds$-almost all  $(s,\omega)\in[t,T]\times\Omega$. From Theorem \ref{Propconv}, Lemma \ref{prop-bound-v}, \reff{XS2} and the monotone convergence theorem, we get by sending $n$ to infinity
\beqs
v(t,x,i) & \geq & \E\Big[\int_{t}^{\nu}f\big(X^{t,x,\alpha}_{s},\alpha_{s}\big)ds-\sum_{t\leq\tau_{k}\leq\nu}c\big(X^{t,x,\alpha}_{\tau_k},\zeta_{k-1},\zeta_{k}\big) +v\big(\nu,X^{t,x,\alpha}_{\nu},\alpha_{\nu}\big)\Big]\;.
\enqs
We now prove the reverse inequality. From 
the definitions of the performance criterion and the value functions, the law of iterated 
conditional expectations and  Markov property of 
our model, we get the successive relations 
\beqs
J(t,x,\alpha) & = & \\
\E\Big[\int_{t}^{\nu}f(s,X^{t,x,\alpha}_{s},\alpha_{s})ds-\sum_{t\leq\tau_{k}\leq\nu}c(X^{t,x,\alpha}_{\tau_k},\zeta_{k-1},\zeta_{k})  & & \\
 +\E^{}\Big[g(X^{t,x,\alpha}_{T})+\int_{\nu}^{T}f(X^{t,x,\alpha}_{s},\alpha_{s})ds-\sum_{\nu<\tau_{k}\leq T}c(X^{t,x,\alpha}_{\tau_k},\zeta_{k-1},\zeta_{k})\Big|\Fc_{\nu}\Big]\Big] & = &\\
\E\Big[\int_{t}^{\nu}f\big(X^{t,x,\alpha}_{s},\alpha_{s}\big)ds-\sum_{t\leq\tau_{k}\leq\nu}c\big(X^{t,x,\alpha}_{\tau_k},\zeta_{k-1},\zeta_{k}\big)+J\big(\nu, X^{t,x,\alpha}_{\nu},\alpha\big)\Big]  & \leq & \\
\E^{}\Big[\int_{t}^{\nu}f\big(X^{t,x,\alpha}_{s},\alpha_{s}\big)ds-\sum_{t\leq\tau_{k}\leq\nu}c\big(X^{t,x,\alpha}_{\tau_k},\zeta_{k-1},\zeta_{k}\big)+v\big(\nu, X^{t,x,\alpha}_{\nu},\alpha_{\nu}\big)\Big]\;.
\enqs
Since $\nu$  and $\alpha$ are arbitrary, we obtain the required inequality.
\ep

\subsection{Viscosity properties}

We prove in this section that the function $v$ is a solution to a system of variational inequalities. More precisely we consider the following PDE
\beq
 \min\Big[-\Dt{v}-\mathcal{L}v-f, v-\Hc v\Big]=0 & \mbox{ on
}& [0,T)\times\Dc\times\Ic,\label{EDPD1}\\
\min\Big[v-g, v-\Hc v\Big]=0 
& \mbox{ on } & \{T\}\times\Dc\times\Ic.\label{EDPDT}
\enq
 where $\Lc$ is
the second order local operator defined by
\beqs
\Lc
v(t,x,i) & =  & \Big(\mu\trans Dv+\frac{1}{2}\mbox{tr}[\sigma\sigma\trans 
D^2v)]\Big)(t,x,i)
\enqs
and $\Hc$ is the nonlocal operator defined by
\beqs
\Hc v(t,x,i) & = & \max_{
\tiny{\begin{array}{c}
j\in\Ic\\
j \neq i
\end{array}}
}\big[v(t,x,j)-c(x,i,j)\big]
\enqs
for all $(t,x,i)\in[0,T]\times\Dc\times\Ic$. As usual, the value functions need not be smooth, and
even not known to be continuous a priori. So, we shall work with the
notion of (discontinuous) viscosity solutions (see \cite{craishlio92}).
Generally, for PDEs arising in optimal control problems involving state constraints, we need  the notion of constrained viscosity solution  introduced by \cite{Son86} for first order equations to take into account the boundary conditions induced by the state constraints.

For a locally bounded function $u$
on $[0,T]\times\Dc \times \Ic$, we define its lower semicontinuous (lsc for
short) envelope $u_*$, and upper semicontinuous (usc for short) envelope $u^*$
by
 \beqs u_*(t,x,i) \; = \; \liminf_{\tiny{\begin{array}{c} (t',x')\rightarrow (t,x),\\(t',x')\in[0,T)\times \Dc \end{array}}} u(t',x',i), & & u^*(t,x,i) \; = \;
\limsup_{\tiny{\begin{array}{c} (t',x')\rightarrow (t,x),\\(t',x')\in[0,T)\times \Dc \end{array}}} u(t',x',i). \enqs
for all $(t,x,i)\in[0,T]\times\Dc\times\Ic$.

\vspace{2mm}

\begin{Remark}\label{lienvv^*}{\rm From Corollary \ref{corvusc} and the definition of the usc envelope we have $v=v^*$ on $[0,T)\times\Dc\times\Ic$. However, this equality may not to be true on $\{T\}\times\Dc\times\Ic$. }
\end{Remark}

\ni We now give the definition of a constrained viscosity solutions to \reff{EDPD1}  and \reff{EDPDT}.

\begin{Definition}
[Constrained viscosity solutions to \reff{EDPD1}-\reff{EDPDT}]

\hspace{1mm}
\begin{enumerate}[(i)]
\item A function $u$, lsc (resp. usc) on $[0,T)\times\Dc\times\Ic$,
is called a viscosity super-solution on $[0,T)\times Int(\Dc)\times\Ic$ (resp. sub-solution on $[0,T)\times\Dc\times\Ic$) to
\reff{EDPD1}-\reff{EDPDT} if 
we have
\beqs
\min\Big[ - \Dt{\varphi}(t,x,i) -   \Lc\varphi(t,x,i)-f(x,i) \; , 
\;  u(t,x,i)-\Hc u(t,x,i) \Big]  & \geq~ (\mbox{resp. } \leq) &  0
\enqs
for any $(t,x,i)\in[0,T)\times Int(\Dc)\times\Ic$ (resp. $(t,x,i)\in[0,T)\times\Dc\times\Ic$), and any $\varphi\in C^{1,2}([0,T]\times\R^d,\R)$ such that
\beqs
\varphi(t,x)-u(t,x,i) & = & \max_{[0,T]\times\Dc} (\varphi-u(.,i)) ~\big(\mbox{resp. }\min_{[0,T]\times\Dc} (\varphi-u(.,i)) \big)
\enqs
and 
\beqs
\min\Big[ u(T,x,i) -g(x,i)~,~  u(T,x,i)-\Hc u(T,x,i) \Big]  & \geq (\mbox{ resp. } \leq) &  0
\enqs
for any $x\in Int(\Dc)$ (resp. $x\in \Dc$).
\item A locally bounded function $u$ on $[0,T]\times\Dc\times\Ic$ is called a constrained viscosity solution to \reff{EDPD1}-\reff{EDPDT} if its lsc envelope $u_*$  is a viscosity super-solution to \reff{EDPD1}-\reff{EDPDT} on $[0,T]\times Int(\Dc)\times\Ic$
and its usc envelope $u^*$ is a viscosity sub-solution  on $[0,T]\times\Dc\times\Ic$ to \reff{EDPD1}-\reff{EDPDT}.
\end{enumerate}
\end{Definition}

\noindent We can now state the viscosity property of $v$.
\begin{Theorem}\label{thvp}
Suppose that the function $v$ is locally bounded. Under \textbf{(H1)} and \textbf{(H2)}, $v$ is a constrained 
viscosity solution to \reff{EDPD1}-\reff{EDPDT}.
\end{Theorem}

\vspace{2mm}

\ni \textbf{Proof of the super-solution property on $[0,T)\times Int(\Dc)\times\Ic$.}
First, for any $(t,x,i)$ $\in$ $[0,T)\times\Dc\times\Ic$, we see, as a consequence of \reff{DPP1} applied to $\nu$ $=$ $t$, and by choosing any admissible control $\alpha$ $\in$ $\Ac^{\Dc}_{t,x,i}$ with immediate switch  $j$ at $t$, that 
\beq\label{vgeqHcv}
v(t,x,i) & \geq & \Hc v(t,x,i)\;.
\enq
Now, let $(\bar t , \bar x,i)\in[0,T)\times Int(\Dc)\times\Ic$ and   $\varphi$ $\in$ $C^{1,2}([0,T]\times \R^d,\R)$ s.t. 
\beq\label{condfctestsursol}
\varphi(\bar t , \bar x)-v_{*}(\bar t , \bar x,i) & = & \max_{[0,T]\times\Dc}(\varphi-v_{*}(.,i)).
\enq
Since $v$ $\geq$ $\Hc v$ on $[0,T)\times Int(\Dc)\times\Ic$, we get from the definition of the operator $\Hc$ and \textbf{(H2)} (i) 
\beqs
v_*(\bar t, \bar x,j ) & \geq & v_*(\bar t, \bar x,j )-c(\bar x,i,j)\;,
\enqs
for all $j\in\Ic$. Therefore we obtain
\beqs
v_{*}(\bar t, \bar x,i ) & \geq & \Hc v_{*}(\bar t, \bar x,i )\;.
\enqs
So it remains to show that
\beq\label{FMSS}
-\Dt{\varphi}(\bar t, \bar x ,i )-\Lc\varphi(\bar t, \bar x, i )-f( \bar x, i ) & \geq & 0\;.
\enq
From the definition of $v_{*}$ there exists a sequence $(t_{m},x_{m})_{m}$ valued in  $[0,T)\times Int(\Dc)$ s.t.  
\beqs
(t_{m},x_{m},v(t_{m},x_{m},i)) & \xrightarrow[m\rightarrow\infty]{} &  (\bar t,\bar x,v_{*}(\bar t, \bar x, i))  \;. 
\enqs
 By continuity of $\varphi$, $\gamma_{m}:= v(t_{m},x_{m},i)-\varphi(t_{m},x_{m})-v_*(\bar t,\bar x,i)+\varphi(\bar t,\bar x)$  converges to 0 as $m$ goes to infinity. Since $(\bar t, \bar x)$ $\in$ $[0,T)\times Int(\Dc)$, there exists $\eta$ $>$ $0$ s.t. for $m$ large enough, $t_{m}$ $<$ $T$ and 
\beqs
((t_{m}-\frac{\eta}{2})\wedge 0,t_{m}+\frac{\eta}{2})\times B(x_{m},\frac{\eta}{2})  & \subset & ((t-\eta)\wedge0,t+\eta)\times B(x,\eta)~~\subset~~ [0,T)\times Int(\Dc)\;.
\enqs
Let us consider an admissible control $\alpha^m$ in $\Ac^{\Dc}_{t_{m},x_{m},i}$ with no switch until the first exit time $\tau_{m}$ before $T$ of the associated process $(s,X^m_{s})$ $:=$ $(s,X^{t_{m},x_{m},\alpha^m}_{s})$ from $(t_{m}-\frac{\eta}{2},t_{m}+\frac{\eta}{2})\times B(x_{m},\frac{\eta}{2})$:
\beqs
\tau_{m} & := & \inf\big\{s\geq t_{m}~:~(s-t_{m})\vee|X^{m}_{s}-x_{m}|\geq \frac{\eta}{2}\big\}\;.
\enqs

Consider also a strictly positive sequence $(h_{m})_{m}$ s.t. $h_{m}$ and $\gamma_{m}/h_{m}$ converge to $0$ as $m$ goes to infinity. By using the dynamic programming principle \reff{DPP1} for $v(t_{m}, x_{m},i)$ and $\nu=\hat \tau _{m}$ $:=$ $\inf\{s\geq t_{m}~:~(s-t_{m})\vee|X^{m}_{s}-x_{m}|\geq \frac{\eta}{4}\}\wedge (t_{m}+h_{m})$, we get 
\beqs
v(t_{m},x_{m},i) & = & \gamma_{m}+v_*(\bar t,\bar x,i)-\varphi(\bar t,\bar x,i) +\varphi(t_{m},x_{m},i) \\
& \geq & \E\Big[\int_{t_{m}}^{\hat \tau _{m}}f(X_{s}^{m},i)ds+ v\big(\hat \tau _{m}, X_{\hat \tau _{m}}^{m},i\big)\Big]\;.
\enqs
Using \reff{condfctestsursol}, we obtain
\beqs
v(t_{m},x_{m},i)  & \geq & \E^{}\Big[\int_{t_{m}}^{\hat \tau _{m}}f( X_{s}^{m},i)ds+\varphi\big(\hat \tau _{m}, X_{\hat \tau _{m}}^{m}\big)\Big]\;.
\enqs
Applying It\^o's formula to $\varphi(s, X_{s}^{m})$ between $t_{m}$ and $\hat \tau _{m}$ and since $\sigma(X_s^m,i)D\varphi(s,X_s^m)$ is bounded for $s\in[t_m,\hat\tau_m]$, we obtain
\beq\label{SS}
\frac{\gamma_{m}}{h_{m}}+ \E^{}\Big[\frac{1}{h_{m}}\int_{t_{m}}^{\hat \tau _{m}}\Big(-\Dt{\varphi}-\Lc \varphi-f\Big)(s, X_{s}^{m},i)ds\Big] & \geq & 0\;,\enq
for all $m\geq 1$. From the continuity of the process $X^{m}$, we have 
\beqs
\P\Big(\exists m,~\forall m'\geq m ~:~\hat \tau _{m'} = t_{m'}+h_{m'}\Big) & = & 1\;.
\enqs 
Hence, by the mean-value theorem, the random variable inside the expectation in  \reff{SS} converges a.s. to $(-\Dt{\varphi}-\Lc \varphi-f)(\bar t,\bar x,i)$ as $m$ goes to infinity. We conclude by the dominated convergence theorem and get \reff{FMSS}.\ep

\vspace{2mm}

\noindent\textbf{Proof of the sub-solution property on $[0,T)\times\Dc\times \Ic$.}  We first recall that $v^*=v$ on $[0,T)\times\Dc\times\Ic$ from Remark \ref{lienvv^*}. 
Let  $(\bar t , \bar x,i) \in [0,T)\times\Dc\times \Ic$ and   $\varphi \in C^{1,2}([0,T]\times\R^d,\R)$ s.t. 
\beq\label{condfctestsousol}
\varphi(\bar t , \bar x)-v(\bar t , \bar x,i) & = & \min_{[0,T]\times\Dc}(\varphi-v(.,i)).
\enq
If $v(\bar t , \bar x,i) \leq \Hc v(\bar t , \bar x,i)$ then the sub-solution property  trivially holds. Consider now the case   $v(\bar t , \bar x,i) > \Hc v(\bar t , \bar x,i)$ and argue by contradiction by assuming on the contrary that 
\beqs
\eta & := & -\Dt{\varphi}(\bar t , \bar x)-\Lc\varphi(\bar t , \bar x,i)-f( \bar x,i) ~ > ~ 0\;.
\enqs
By continuity of $\varphi$ and its derivatives, there exists some $\delta$ $>$ $0$ such that $\bar t + \delta$ $<$ $T$ and 
\beq
\Big(-\Dt{\varphi}-\Lc\varphi-f\Big)( t ,  x,i)  & \geq & \frac{\eta}{2}\;, 
\enq
for all $(t,x)\in\Vc:=\Big((\bar t -\delta,\bar t +\delta)\cap[0,T)\Big)\times B(\bar x,\delta)$.
 By the dynamic programming principle \reff{DPP1}, given $m$ $\geq$ $1$, there exists $\hat \alpha ^m$ $=$ $(\hat \tau ^m_{n},\hat \zeta ^m_{n})_{n}\in\Ac^{\Dc}_{\bar t,\bar x,i}$ s.t. for any stopping time $\tau$ valued in $[\bar t,T]$, we have 
\beqs
v(\bar t,\bar x,i) & \leq & \E^{}\Big[\int_{\bar  t}^{\tau}f(\hat X ^m_{s},i)-\sum_{\bar  t\leq\hat \tau^m_{n}\leq \tau}c(\hat X ^m_{\hat \tau^m_{n}},\hat \zeta _{n}^m,\hat \zeta _{n}^m)+v(\tau,\hat X ^m_{\tau},i)\Big]+\frac{1}{m}
\enqs
where $\hat X ^m:=X^{\bar t,\bar x,\hat \alpha ^m}$.  By choosing $\tau$ $=$ $\bar \tau _{m}$ $:=$ $\hat \tau ^m_{1} \wedge \nu ^m$ where
\beqs
\nu ^m & := & \inf\{ s\geq \bar t ~:~(s,\hat X ^m_{s})\notin \Vc \}
\enqs 
is the first exit time of $(s,\hat X ^m_s)$ from $\Vc$, 
 we then get
\beq
v(\bar t,\bar x,i)
& \leq & \E^{}\left[\int_{\bar t}^{\bar \tau ^m} f(\hat X^m_{s},i)ds\right]+\E^{}\left[ v(\bar \tau ^m,\hat X ^m_{\bar \tau ^{m}},i  )\mathds{1}_{\nu ^m  < \hat \tau ^m_{1} } \right]\nonumber\\
 & & + \E^{}\left[ [v(\bar \tau ^m,\hat X ^m_{\bar \tau ^{m}},\hat \zeta ^m_{1})-c(\hat X ^m_{\bar \tau ^{m}},i,\hat \zeta ^m_{1})]\mathds{1}_{\nu ^m  \geq \hat \tau ^m_{1} } \right]+\frac{1}{m}\nonumber\\
& \leq &\nonumber \E^{}\left[\int_{\bar t}^{\bar \tau ^m} f(\hat X^m_{s},i)ds\right]+\E^{}\left[ v(\bar \tau ^m,\hat X ^m_{\bar \tau ^{m}} ,i )\mathds{1}_{\nu ^m  < \hat \tau ^m_{1} } 
\right]\\
 & & \label{enegvppdsssol}+ \E^{}\left[ \Hc v(\bar \tau ^m,\hat X ^m_{\bar \tau ^{m}},i  )\mathds{1}_{\nu ^m  \geq \hat \tau ^m_{1} } \right]+\frac{1}{m}\;.\quad\qquad\label{inegv}
\enq
Now, since $v$ $\geq$ $\Hc v$ on $[0,T]\times\Dc\times\Ic$ and $\hat \alpha^m\in\Ac^{\Dc}_{\bar t,\bar x,i}$, we obtain from \reff{condfctestsousol}
\beqs
\varphi(\bar t,\bar x,i)  & \leq & \E^{}\left[\int_{\bar t}^{\bar\tau_{m}}f(\hat X ^m_{s},i)ds+ \varphi(\bar \tau ^m,\hat X ^m_{\bar \tau ^{m}}) \right]+\frac{1}{m}\;.
\enqs 
Applying It\^o's formula to $\varphi(s,\hat X ^m_{s})$ between $t_{m}$ and $\bar \tau ^m$ we get:
\beqs
0 & \leq &  \E^{}\left[ \int_{t_{m}}^{\bar \tau ^m}(\Dt{\varphi}+\Lc\varphi+f)(s,\hat X ^m_{s},i) \right]+\frac{1}{m}
~\leq~-\frac{\eta}{2}\E^{}\big[\bar \tau ^{m}-\bar t\big]+\frac{1}{m}\;.
\enqs
This implies 
\beq\label{cvgtaum}
\lim_{m\rightarrow+\infty} \E^{}[\bar \tau ^{m}] & =  & \bar t\;.
\enq
From the definition of $\nu^m$ and \reff{cvgtaum} we have, up to a subsequence,  
\beq\label{conv-proba-taudelata-geq-tau1}
\P\big(\nu ^m  \geq \hat \tau ^m_{1} \big) & \xrightarrow[m\rightarrow\infty]{}  & 1\;.
\enq
On the other hand, we get from \reff{enegvppdsssol}
\beqs
v(\bar t,\bar x,i) & \leq &\E^{}\left[\int_{\bar  t}^{\bar\tau_{m}}f(\hat X ^m_{s},i)ds\right] + \P\big(\nu ^m  < \hat \tau ^m_{1} \big)\sup_{(t',x')\in Adh(\Vc)}v(t',x',i)\\ 
 &  & 
+ \;\P\big(\nu ^m  \geq \hat \tau ^m_{1} \big)\sup_{(t',x')\in Adh(\Vc)}\Hc v(t',x',i)+\frac{1}{m}\;.
\enqs
From Lemma \ref{prop-bound-v}, \reff{cvgtaum} and \reff{conv-proba-taudelata-geq-tau1} we get by sending $m$ to $\infty$
\beqs
v(\bar t,\bar x,i) & \leq & \sup_{(t',x')\in Adh(\Vc)}\Hc v(t',x',i)\;.
\enqs
Since $v=v^*$, we get by sending $m$ to infinity and $\delta$ to zero 
\beqs
v(\bar t,\bar x,i) & \leq & (\Hc v)^*(\bar t,\bar x,i)~\leq~ \Hc v(\bar t,\bar x,i)\;,
\enqs 
which is the required contradiction.\ep

\vspace{2mm}

\ni\textbf{Proof of the viscosity super-solution property on $\{T\}\times Int(\Dc)\times\Ic$.} 
Fix some $(\bar x,i)$ $\in$ $Int(\Dc)\times\Ic$, and consider a sequence 
$(t_{m},x_{m})_{m\geq1}$ valued in  $[0,T)\times Int(\Dc)$, such that
\beqs
\big(t_m,x_m,v(t_m,x_m,i)\big) & \xrightarrow[m\rightarrow\infty]{} & (T,\bar x,v_*(T,\bar x,i)\big)\;. 
\enqs 
Let $\delta>0$ s.t. $B(\bar x,\delta)\in Int(\Dc)$. We first can suppose w.l.o.g. that
\beq\label{emboite}
B(x_m,{\delta\over 2}) & \subset & B (\bar x, \delta)
\enq
for all $m\geq 1$. By taking a strategy $\alpha^m=(\tau_k^{m},\zeta_k^{m})_{k}\in\Ac_{t_{m},x_{m},i}^{\Dc}$ with no switch before $\nu_{m}$ $:=$ $\inf\{s\geq t_{m}, ~X_{s}^{m} \notin B(x_m,{\delta\over 2})\}\wedge T$ with $X^{m}$ $:=$ $X^{t_{m},x_{m},\alpha^m}$, we have from \reff{DPP1} applied to $\tau_{m}$ $:=$ $\inf\{s\geq t_{m}, ~X_{s}^{m} \notin B(x_m,{\delta\over 4})\}\wedge T$ and $\alpha_{m}$  
\beqs
v(t_m,x_m,i) & \geq & \E^{}\Big[ \int_{t_{m}}^{\tau^m}f(X^m_{s},i)ds\Big]+\E^{}\big[v(\tau^m,X^m_{\tau^{m}},i)\big]
\enqs
Since $v(T,.)=g$ we obtain from \reff{emboite}
\beq\nonumber
v(t_m,x_m,i) & \geq & \E^{}\Big[ \int_{t_{m}}^{\tau^m}f(X^m_{s},i)ds\Big]+\E^{}\Big[v(\tau^m,X^m_{\tau^{m}},i)\mathds{1}_{\tau^m<T}\Big]+\E^{}\Big[g(X^m_{\tau^{m}},i)\mathds{1}_{\tau^m=T}\Big]\\\nonumber
 & \geq & \E^{}\Big[ \int_{t_{m}}^{\tau^m}f(X^m_{s},i)ds\Big] + \P\big(\tau^m<T\big)\inf_{\tiny{\begin{array}{c}t<T\\
 x\in Adh(B(\bar x,{\delta}))
 \end{array}}}v(t,x,i) \\
 & & +\P\big(\tau^m=T\big)\inf_{x\in Adh( B(\bar x,{\delta}))}g(x)\;.\label{majv>g}
\enq
Since $\E^{}[\sup_{s\in [t_{m},T]}|X^m_{s}-x_{m}|]$  converges to zero (see e.g. Corollary 12, Section 5, Chapter 2 in \cite{Kryl80}), we have, up to a subsequence, 
\beqs
\sup_{s\in [t_{m},T]}|X^m_{s}-x_{m}| &  \xrightarrow[m\rightarrow\infty]{\P-a.s.}   & 0\;.
\enqs
From the convergence of  $(x_{m})_{m}$ to $x\in Int(\Dc)$, we deduce that 
\beqs
\P\big(\tau^m=T\big) & \xrightarrow[m\rightarrow\infty]{} & 1 \;.
\enqs
Sending $m$ to infinity and $\delta$ to $0$ in \reff{majv>g} we get
\beq \label{vgeqU}
v_*(T,\bar x,i) & \geq & g(\bar x,i)\,. 
\enq
On the other hand, we know from \reff{vgeqHcv} that $v$ $\geq$ $\Hc v$ on $[0,T)\times Int(\Dc)$, and thus 
\beqs
v(t_m,x_m,i) &\geq&  \Hc v(t_m,x_m,i) \; \geq \; \Hc v_*(t_m,x_m,i), 
\enqs 
for all  $m \geq 1$. Recalling that $\Hc v_*$ is lsc, we obtain by sending $m$ to infinity 
\beqs
v_*(T,\bar x,i) & \geq & \Hc v_*(T,\bar x,i). 
\enqs
Together with \reff{vgeqU}, this proves the required viscosity super-solution property of \reff{EDPDT}. \ep

\vspace{2mm}

\ni\textbf{Proof of the viscosity sub-solution property on $\{T\}\times\Dc\times\Ic$.} 
We 
 argue by contradiction by assuming that there exists $(\bar x,i)$ $\in$ $\Dc \times\Ic $ such that 
\beq \label{vsous}
\min\big[ v^*(T,\bar x,i)  - g(\bar x,i)~,\Hc v^*(T,\bar x,i) \big] &:=& 2\eps~>~0.
\enq
One can find a sequence of smooth functions $(\varphi^n)_{n\geq 0}$ on $[0,T]\times\R^d$ such that $\varphi^n$ converges pointwisely to $v^*(.,i)$ on $[0,T]\times\Dc\times\Ic$  as $n\rightarrow\infty$. Moreover, by \reff{vsous} and the upper semicontinuity of $v^*$, 
we may assume that the inequality 
\beq\label{maj1}
\min\big[ \varphi^n  - g(.,i)~,~\varphi^n - \max_{j\in\Ic} \{v^*(.,j)+c(.,i,j)\}  \big] & \geq & \eps,
\enq
holds on some bounded neighborhood $B^{n}$ of $(T,\bar x)$ in $[0,T]\times\Dc$, for $n$ large enough.  
Let  $(t_{m},x_{m})_{m\geq1}$ be a sequence in $[0,T)\times \Dc$ such that 
\beqs
\big(t_m,x_m,v(t_m,x_m,i)\big) & \xrightarrow[m\rightarrow\infty]{} & (T,\bar x,v^*(T,\bar x,i)\big)\;. 
\enqs
Then there exists $\delta^n$ $>$ $0$ such that  
$B^{n}_m$ $:=$ $[t_{m},T]\times B(x_{m},\delta^{n})$ $\subset$ $B^n$ for $m$ large enough, so that \reff{maj1} holds on  $B^{n}_m$.  Since $v$ is locally bounded, there exists  some $\eta>0$ such that $|v^*|$ $\leq$ $\eta$ on $B^{n}$. We can then assume that $\varphi^n$ $\geq$ $-2\eta$ on $B^{n}$. Let us define the smooth function $\tilde \varphi _{m}^n$ by
\beqs
\tilde \varphi _{m}^n(t,x) & := &  \varphi ^n(t,x)+ \Big(4\eta\frac{|x-x_{m}|^2}{|\delta^{n}|^2}+ \sqrt{T-t}\Big)
\enqs
for $(t,x)\in[0,T)\times Int(\Dc)$ 
and observe that 
\beq\label{maj2}
(v^*-\tilde \varphi _{m}^n)(t,x,i) & \leq & -\eta, 
\enq
for   $ (t,x)\in[t_{m},T]\times\partial B(x_{m},\delta^{n})$. 
Since $\Dt{\sqrt{T-t}}\longrightarrow-\infty$ as $t\rightarrow T$, we have for $m$ large enough
\beq\label{maj3}
-\Dt{\tilde \varphi _{m}^n}  - \Lc \tilde \varphi _{m}^n(.,i) & \geq & 0 ~ \mbox{ on }~ B^{n}_m. 
\enq
Let $\alpha^m=(\tau^m_{j},\zeta^m_{j})_{j}$ be a $\frac{1}{m}-$optimal control for $v(t_{m},x_{m},i)$ with corresponding state process 
$X^m=X^{t_m,x_m,\alpha^m}$, and denote by   $\theta_n^m=\inf\big\{ s\geq t_{m}~:~ (s,X_{s}^m)\notin B^{n}_m \big\}\wedge \tau_{1}^m\wedge T$. 
From \reff{DPP1} we have
\beq
v(t_{m},x_{m},i) - \frac{1}{m} & \leq & \E^{}\Big[\int_{t_{m}}^{\theta_{n}^m}f(X_{s}^m,i)ds\Big] 
+  \E^{}\Big[\mathds{1}_{\theta_n^m<\tau_{1}^m\wedge T} \; 
v(\theta_n^m,X^m_{\theta_n^m},i)\Big] \label{DPPinter} \\
& & \; + \;  
\E^{}\Big[\mathds{1}_{\theta_n^m=T < \tau_{1}^m} \;  g(X^m_{\theta_n^m},i)\Big]\nonumber\\
& &  
+ \E^{}\Big[\mathds{1}_{\tau_{1}^m = \theta_n^m\leq T} \; \Big(v\big(\tau_{1}^m,X_{\tau_{1}^m}^m,\zeta^m_{1}\big)+c(X_{\tau_{1}^m}^m,i,\zeta^m_{1})\Big)\Big]\;.
\nonumber
\enq
Now,  by applying It\^{o}'s Lemma to $\tilde \varphi _{n}^m(s,X_{s}^m)$ between $t_{m}$ and $\theta_n^m$
we get   from  \reff{maj1}, \reff{maj2} and \reff{maj3}
\beqs \label{inegphi}
\tilde \varphi^n_{m}(t_{m},x_{m}) & \geq & 
\E^{}\Big[ \mathds{1}_{\theta_n^m<\tau_{1}^m} \; \tilde \varphi^n_{m}(\theta_n^m,X^m_{\theta_n^m}) \Big] 
+ \E^{}\Big[ \mathds{1}_{\tau_{1}^m\leq \theta_n^m}\tilde \varphi^n_{m}\big(\tau_{1}^m,X^m_{\tau_{1}^m}\big) \Big]\nonumber\\
 & \geq & \E^{}\Big[\mathds{1}_{\theta_n^m<\tau_{1}^m\wedge T}\Big( v^*(\theta_n^m,X^m_{\theta_n^m},i)+\eta \Big)\Big]
 +\E^{}\Big[\mathds{1}_{\theta_n^m=T < \tau_{1}^m}\Big(g(X^m_{\theta_n^m},i)+\eps \Big)\Big]\nonumber\\
 &  & \; + \;  \E^{}\Big[\mathds{1}_{\tau_{1}^m = \theta_n^m\leq T} 
 \Big(v^*\big(\tau_{1}^m,X_{\tau_{1}^m}^m,\zeta^m_{1}\big)+c(X_{\tau_{1}^m}^k,i,\zeta_{1}^m)  +\eps \Big) \Big]. 
\enqs
Together with \reff{DPPinter}, this implies
\beqs
\tilde \varphi^n_{m}(t_{m},x_{m}) & \geq & v(t_{m},x_{m},i) - \E^{}\Big[\int_{t_{m}}^{\theta_{n}^m}f(X_{s}^m,i)ds\Big] 
 -\frac{1}{m} +\eps\wedge\eta.
\enqs
Sending $m$, and then $n$ to infinity,  we get the required  contradiction: $v^*(T,\bar x,i)$ $\geq$ $v^*(T,\bar x,i) + \eps\wedge\eta$.  
\ep

\section{Uniqueness result}

\subsection{ Maximality of the value function as a solution to the SVI}
In general, the uniqueness of a viscosity solution to some PDE is given by a comparison theorem. Such a result says that for $u$ an $usc$ super-solution and and $w$ a $lsc$ sub-solution, we have $u\geq w$. Applying this result to $u=v_*$ the $lsc$ envelope of $v$ and $w=v^*$ the $usc$ envelope of $v$ we would get that $v_* = v^*$ and $v$ would be continuous. As the counter-example presented in Subsection \ref{contreexemple} shows, such a property cannot hold for SVI \reff{EDPD1}-\reff{EDPDT}. 

\vspace{2mm}

We therefore provide a weaker characterization of $v$. 
To this end, we introduce, for $n\geq 1$, the SVI with penalized coefficients defined on the whole space $[0,T]\times\R^d\times\Ic$:
\beq
 \min\Big[-\Dt{v}-\mathcal{L}v-f_n, v-\Hc v\Big]=0 & \mbox{ on
}& [0,T)\times\R^d\times\Ic,\label{EDPDWD1}\\
\min\Big[v-g_n, v-\Hc v\Big]=0 
& \mbox{ on } & \{T\}\times \R^d\times\Ic.\label{EDPDWDT}
\enq
Under assumption \textbf{(H1)} and \textbf{(H2)}, we can use Lemma \ref{Propfngn} to apply Proposition 5.1 in \cite{bou06} and we get from Proposition 4.12 in \cite{bou06} the following  comparison result for this PDE.

\begin{Theorem}\label{compthRd} Suppose that \textbf{(H1)} and \textbf{(H2)} hold.
Let $u$ and $w$ be respectively a sub-solution and a super-solution to \reff{EDPDWD1}-\reff{EDPDWDT}. Suppose that there exists two constants $C_u>0$ and $C_w>0$ and an integer $\gamma\geq1$ such that
\beqs
u(t,x,i) & \leq & C_u\big(1+|x|^\gamma\big)\\
w(t,x,i) & \geq & -C_w\big(1+|x|^\gamma\big)
\enqs
for all $(t,x,i)\in[0,T]\times\R^d\times \Ic$. Then 
we have $u\leq w$ on $[0,T]\times\R^d\times \Ic$.
\end{Theorem}

\ni We now introduce the following additional assumption on the function $v$.

\vspace{2mm}

\ni \textbf{(H3)} There exists a constant $C>0$ and an integer $q\geq 1$ such that
\beq\label{condCrLin}
v(t,x,i) & \geq & -C\big(1+|x|^q\big) 
\enq
for all $(t,x,i)\in[0,T]\times\Dc\times \Ic$.

\vspace{2mm}

\ni We give in the next subsection, some examples where \textbf{(H3)} is satisfied. 
 We can state our maximality result as follows.
\begin{Theorem}
Under \textbf{(H1)}, \textbf{(H2)} and \textbf{(H3)} the function $v$ is the maximal constrained viscosity solution to \reff{EDPD1}-\reff{EDPDT} satisfying \reff{condCrLin}: for any function $w:[0,T]\times\Dc\times\Ic\rightarrow\R$ such that
\begin{itemize}
\item $w$ is a constrained viscosity solution to \reff{EDPD1}-\reff{EDPDT}, 
\item there exists a constant $C$ and an integer $\eta\geq1$ such that 
\beq\label{condCrLinw}
w(t,x,i) & \geq & -C\big(1+|x|^\eta\big)
\enq
for all $(t,x,i)\in[0,T]\times\Dc\times\Ic$,
\end{itemize}
 we have $v\geq w$ on $[0,T]\times\Dc\times\Ic$.
\end{Theorem}
\ni \textbf{Proof.} Let $w:[0,T]\times\Dc\times\Ic\rightarrow\R$ be a constrained viscosity solution to \reff{EDPD1}-\reff{EDPDT} satisfying \reff{condCrLinw}. We proceed in four steps to prove that $w\leq v$.

\vspace{2mm}

\ni \textbf{Step 1.} \textit{Extension of the definition of $w$ to $[0,T]\times\R^d\times\Ic$}.

\ni  For $n\geq 1$, we define the function $\tilde w_n$ on $[0,T]\times\R^d\times \Ic$ by
\begin{equation}\label{ext-w}
\tilde w_n(t,x,i)  =  \left\{
\begin{array}{ccl}
w(t,x,i)  & \mbox{ for } &  (t,x,i)\in[0,T]\times\Dc\times\Ic\;,\\
-C_ne^{-\rho_n t}\big(1+|x|^{2\eta}\big) & \mbox{ for } &  (t,x,i)\in[0,T]\times(\R^d\setminus\Dc)\times\Ic\;.
\end{array}
\right. 
\end{equation}
where $\rho_n$ and $C_n$ are two positive constants.
From \textbf{(H1)}, \textbf{(H2)}, Lemma \ref{Propfngn} and \reff{condCrLinw}, we can find $\rho_n$ and $C_n$ (large enough) such that
\beq\label{cond1tildew}
-\Dt{\tilde w_n}-\mathcal{L}\tilde w_n-f_n &\leq & 0~ \mbox{ on
}~[0,T)\times(\R^d\setminus\Dc)\times\Ic\;,\\
\label{cond2tildew}
\tilde w_n-g_n & \leq & 0 
~ \mbox{ on } ~ \{T\}\times \R^d\times\Ic\;,
\enq
and 
\beq\label{cond3tildew}
\tilde w_n(t,x,i) & \geq & -C_ne^{-\rho_n t}\big(1+|x|^{2\eta}\big) ~ \mbox{ for
}~(t,x,i)\in [0,T)\times\R^d\times\Ic\;.
\enq


\ni \textbf{Step 2.} \textit{Viscosity property of $\tilde w_n$.}

\ni For $C_n$ and $\rho_n$ such that \reff{cond1tildew},\reff{cond2tildew} and \reff{cond3tildew} hold, we obtain that $\tilde w_n$  is a viscosity sub-solution to \reff{EDPDWD1}-\reff{EDPDWDT}. Indeed, let $\varphi\in C^{1,2}([0,T]\times\R^d,\R)$ and $(t,x,i)\in[0,T]\times\R^d\times\Ic$ such that
\beq\label{condfcttest}
(\tilde w_n^*-\varphi)(t,x,i) & = & \max_{[0,T]\times\R^d\times\Ic}(\tilde w_n^*-\varphi)\;.
\enq
We first notice from \reff{cond3tildew} that the upper semicontinuous envelope $\tilde w_n^*$ of $\tilde w_n$ is given by
\begin{equation}\label{characw^*}
\tilde w_n^*(t,x,i)  =  \left\{
\begin{array}{ccl}
w^*(t,x,i)  & \mbox{ for } &  (t,x,i)\in[0,T]\times\Dc\times\Ic\;,\\
-C_ne^{-\rho_n t}\big(1+|x|^{2\eta}\big) & \mbox{ for } &  (t,x,i)\in[0,T]\times(\R^d\setminus\Dc)\times\Ic\;.
\end{array}
\right. 
\end{equation}
We now prove that $\tilde w_n$ is a sub-solution to  \reff{EDPDWD1}-\reff{EDPDWDT}. Using \reff{cond2tildew}, \reff{characw^*} and the viscosity sub-solution property of $w$, we get
\beqs
\tilde w_n^* & \leq & g_n~~\mbox{ on }~ \{T\}\times\R^d\times\Ic\;.
\enqs
For the viscosity property on $[0,T)\times\R^d\times\Ic$, we distinguish two cases.
\begin{itemize}
\item Case 1: $(t,x,i)\in[0,T)\times\Dc\times\Ic$. From \reff{condfcttest} and \reff{characw^*}, we have
\beqs
(\tilde w_n^*-\varphi)(t,x,i) & = & \max_{[0,T]\times\Dc\times\Ic}(\tilde w_n^*-\varphi)\;.
\enqs
Since $w$ is a constrained viscosity solution to \reff{EDPD1}-\reff{EDPDT} and $f=f_n$ on $\Dc$ we get
\beqs
 \min\Big[-\Dt{\varphi}(t,x,i)-\mathcal{L}\varphi(t,x,i)-f_n(t,x,i), \varphi(t,x,i)-\Hc \tilde w^*_n(t,x,i)\Big] & \leq & 0 \;. 
\enqs
\item Case 2: $(t,x,i)\in[0,T)\times(\R^d\setminus \Dc)\times\Ic$. From \reff{cond1tildew}, \reff{characw^*} we  also get
\beqs
 \min\Big[-\Dt{\varphi}(t,x,i)-\mathcal{L}\varphi(t,x,i)-f_n(t,x,i), \varphi(t,x,i)-\Hc \tilde w_n^*(t,x,i)\Big] & \leq & 0 \;. 
\enqs
\end{itemize}
Therefore, $\tilde w_n$ is a viscosity sub-solution to \reff{EDPDWD1}-\reff{EDPDWDT}.

\vspace{2mm}

\ni \textbf{Step 3.} \textit{Growth condition on $v_n$.}


\ni We prove that for each $n\geq1$ there exists a constant $C_n>0$ such that
\beqs
 v_n(t,x,i) & \geq &  -C_n\big(1+|x|^{2q}\big)\;,\quad (t,x,i)\in[0,T]\times\R^d\times\Ic\;.
\enqs
Fix $(t,x,i)\in[0,T]\times\R^d\times\Ic$, and denote by $^0\alpha=(^0\tau_k,{^0\zeta_k})_k$ the trivial strategy of $\Ac_{t,i}$ $i.e.$ $^0\tau_0=t$, $^0\zeta_0=i$ and $^0\tau_k>T$ for $k\geq 1$. Then we have
\beqs
v_n(t,x,i) & \geq & J_n(t,x,{{}^0\alpha}) 
\enqs
From the definition of $J_n$, \reff{XS2} and Lemma \ref{Propfngn} there exists a constant $\tilde C_n>0$ such that
\beqs
v_n(t,x,i) & \geq & -\tilde C_n\big(1+|x|^q\big)\;. 
\enqs
for all $(t,x,i)\in[0,T]\times\R^d\times\Ic$.

\vspace{2mm}

\ni \textbf{Step 4.} \textit{Comparison on $[0,T]\times\R^d\times\Ic$.}
From Proposition 4.2 in \cite{bou06}, we know that $v_n$ is a viscosity solution to \reff{EDPDWD1}-\reff{EDPDWDT}.
Using the results of Steps 2 and 3, 
we can apply Theorem \ref{compthRd} to $\tilde w_n$ and $v_n$ with $\gamma =2\eta+q$, and we get 
\beqs
\tilde w_n (t,x,i)~~\leq~~\tilde w_n^*(t,x,i) & \leq & v_n(t,x,i)\;,
\enqs
for all $(t,x,i)\in[0,T]\times\R^d\times\Ic$. Sending $n$ to infinity and using Theorem \ref{Propconv} and \reff{ext-w}, we get $w\leq v$ on $[0,T]\times\Dc\times\Ic$.
\ep
\begin{Remark}{\rm We notice that the counter-example given in the  Sub-section \ref{contreexemple} also satisfies Assumption \textbf{(H3)}. In particular this gives an example where the classical uniqueness does not hold and where our maximality result is valid. }
\end{Remark}
\vspace{2mm}

\subsection{Sufficient conditions for \textbf{(H3)}}
We end this Section by providing explicit examples where \textbf{(H3)} is satisfied. The idea consists in constructing switching strategies with finite number of switches and satisfying the constraint imposed on the controlled diffusion. This allows to get a lower bound for the value function. Thanks to the estimate of Lemma \ref{prop-bound-v}, this proves the polynomial growth of the value function.

\vspace{2mm}

 The first example deals with the case where there exists a regime that stops the controlled diffusion. By switching immediately on it, we keep the controlled diffusion stays in $\Dc$. The second example considers the case where for any initial condition there exists an associated regime that keeps the associated diffusion in $\Dc$. By switching on such a regime at the first time the diffusion meets the boundary $\partial \Dc$ of $\Dc,$ we get a strategy satisfying the constraint. Finally, the last example concerns the case of a convex domain $\Dc$.  Using a viability condition involving the normal cone we also ensure the existence of a regime keeping the diffusion in $\Dc$.
We notice that all the presented conditions are satisfied by the examples presented in Section \ref{example elec}.  
\begin{Proposition}
(i) Suppose that for any $x\in \partial \Dc$ there exists $i_x\in\Ic$ such that $\mu(x,i_x)=0$ and $\sigma(x,i_x)=0$, then  \textbf{(H3)} holds.

\vspace{1mm}

\ni(ii) Suppose that for each $(t,x)\in[0,T]\times\Dc$, there exists $i_{t,x}\in\Ic$ such that the process $X^{t,x}$ defined by
\beqs
X_s^{t,x} & = & x+\int_t^s\mu(X_r^{t,x},i_{t,x})dr+\int_t^s\sigma(X_r^{t,x},i_{t,x})dW_r\;,\quad s\geq t\;,
\enqs 
satisifies
\beq\label{condindiceadmissible}
\P\big(X_s^{t,x}\in\Dc,~\forall s\in[t,T]\big) & = & 1\;.
\enq
Then \textbf{(H3)} is satisfied.

\vspace{1mm}

\ni(iii) Suppose that $\Dc$ is convex and  there exists $i^*\in \Ic$ such that
\beqs
p\trans\mu(x,i^*)+ \frac{1}{2}\textrm{tr}[\sigma(x,i^*)\sigma(x,i^*)\trans 
A)] & \leq & 0
\enqs
for all $x\in \partial \Dc$ and all $(p,A)\in\Nc^2_{\Dc}(x)$ where $\Nc^2_{\Dc}(x)$ is the second order normal cone to $\Dc$ at $x$ defined by
\beqs
\Nc_{\Dc}(x) & = & \Big\{ (p,A)\in\R^d\times\S^d~:~p\trans (y-x)+\frac{1}{2}(y-x)\trans A(y-x)\;\leq\;o(|y-x|^2) \\
 & & \qquad\qquad\qquad\qquad\qquad\qquad\qquad\qquad\qquad\qquad \mbox{ as } ~y\rightarrow x ~\mbox{ and } y\in \Dc\Big\}\;\;,
\enqs
and $\S^d$ is the set of $d\times d$ symmetric matrices.. Then  \textbf{(H3)} holds.

\end{Proposition}
\ni \textbf{Proof.} (i) Fix an initial condition $(t,x,i)\in [0,T]\times\Dc\times\Ic$. Let $X^{t,x}$ be the diffusion defined by
\beqs
X^{t,x}_s & = & x+\int_t^s\mu(X^{t,x}_r,i)dr+\int_t^s\sigma(X^{t,x}_r,i)dW_r\;,\qquad s\geq t\;.
\enqs
 Consider the strategy $\alpha:(\tau_k,\zeta_k)_k$ defined by $(\tau_0,\zeta_0) =  (t,i)$,
\beqs
\tau_1 & = & \inf\big\{s\geq 0~:~X_s\in\partial\Dc\big\}\\
\zeta_1 & = & i_{X_{\tau_1}}
\enqs
and $\tau_k>T$ and $\zeta_k=\zeta_1$ for $k\geq 2$. We then have $\mu(X^{t,x,\alpha}_s,\alpha_s)$=0 and $\sigma(X^{t,x,\alpha}_s,\alpha_s)=0$ for $s\in [\tau_1,T]$. Therefore, we get $\alpha\in\Ac_{t,x,i}^{\Dc}$ and
\beqs
v(t,x,i) & \geq & J(t,x,\alpha)\;.
\enqs
From \reff{XS2} and \textbf{(H2)} (ii) there exists a constant $C>0$ such that
\beqs
v(t,x,i) & \geq & -C(1+|x|^q)\;.
\enqs
By combining this inequality with Lemma \ref{prop-bound-v}, we get \textbf{(H3)}.

\vspace{2mm}

\ni (ii) Fix $(t,x,i)\in[0,T]\times\Dc\times\Ic$. Consider the strategy $\alpha=(\tau_k,\zeta_k)_k$ defined by $(\tau_0,\zeta_0) =  (t,i)$,
$(\tau_1,\zeta_1)  =  (t,i_{t,x})$
and $\tau_k>T$ for $k\geq 2$. From \reff{condindiceadmissible} we get $\alpha\in\Ac_{t,x,i}^{\Dc}$.  We then have
\beqs
v(t,x,i) & \geq & J(t,x,\alpha)\;.
\enqs
From \reff{XS2} and \textbf{(H2)} (ii) there exists a constant $C>0$ such that
\beqs
v(t,x,i) & \geq & -C(1+|x|^q)\;.
\enqs
This inequality with Lemma \ref{prop-bound-v} give \textbf{(H3)}.

\vspace{2mm}

\ni (iii) From Proposition 8 and Remark 9 in \cite{GIS13} we get that for any initial condition $(t,x,i)\in [0,T]\times\Dc\times\Ic$, the control $\alpha=(\tau_k,\zeta_k)_k$ defined by
\beqs
(\tau_0,\zeta_0) & = & (t,i)\\
(\tau_1,\zeta_1) & = & (t,i^*)
\enqs  
and $\tau_k>T$ for $k\geq2$, satisfies $\alpha\in\Ac_{t,x,i}^{\Dc}$.  We then have
\beqs
v(t,x,i) & \geq & J(t,x,\alpha)\;.
\enqs
From \reff{XS2} and \textbf{(H2)} (ii) there exists a constant $C>0$ such that
\beqs
v(t,x,i) & \geq & -C(1+|x|^q)\;.
\enqs
Using Lemma \ref{prop-bound-v}, we get \textbf{(H3)} from this last inequality  .
\ep

\appendix
\section{Appendix}

\subsection{Additional results on convergence and measurability}\label{apptd. means cvgce}
We first present two results about stopping times and measurability.
 


\begin{Proposition}\label{ecriture-ta} Let $(\Omega,\Gc,\P)$ be a complete probability space endowed with a Brownian motion $B$. 
Let $\mathbb{H}=(\Hc)_{t\geq0}$ be the complete right-continuous filtration generated by $B$, $\tau$ an $\mathbb{H}$-stopping time and $\zeta$ an $\Hc_\tau$-measurable random variable. Suppose that there exists a constant $M$ such that $\P(\tau\leq M)=1$. Then there exist 
two Borel functions $\psi$ and $ \phi$ such that  
\beqs
\tau ~ = ~  \psi\big((B_s)_{s\in[0,M]}\big) & \mbox{ and } & \zeta~=~\phi\big((B_s)_{s\in[0,M+1]}\big) ~~\quad \P-a.s.
\enqs
\end{Proposition}
\noindent\textbf{Proof.} Since $\tau\leq M$ $\P$-a.s. we can write
\beq\label{decompTauIntegrale}
\tau & = & \int_0^M\mathds{1}_{\tau>s}ds~~=~~\lim_{n\rightarrow\infty}{M\over n}\sum_{k=0}^{n-1}\mathds{1}_{\tau>{k\over n}M}\;,\quad \P-a.s.
\enq
Since $\tau$ is a $\mathbb{H}$-stopping time and $\mathbb{H}$ is the complete right-continuous extension of the natural filtration of $B$, we can write from Remark $32$, Chapter 2 in \cite{delmey75}  
\beq\label{ineg indic}
\underline \psi_n^k\big((B_s)_{s\in[0,M]}\big) ~\leq ~\mathds{1}_{\tau> {k\over n}M}~\leq~\bar \psi_n^k\big((B_s)_{s\in[0,M]}\big)
\enq
and 
\beq\label{proba-diff}
\P\Big(\underline \psi_n^k\big((B_s)_{s\in[0,M]}\big)~\neq~ \bar \psi_n^k\big((B_s)_{s\in[0,M]}\big)\Big) & = & 0
\enq
where $\underline \psi_n^k$ and $\bar \psi_n^k$ are two Borel  functions for any $n\geq 1$ and any $k\in\{0,\ldots,n-1\}$. Define the Borel functions $\bar \psi_n$ and $\underline \psi_n$ by  
\beqs
\bar \psi_n ~~ = ~~ {M\over n}\sum_{k=0}^{n-1} \bar \psi_n^k  & \mbox{ and } & \underline \psi_n ~~ = ~~ {M\over n}\sum_{k=0}^{n-1} \underline \psi_n^k 
\enqs
We then get from \reff{decompTauIntegrale}, \reff{ineg indic} and \reff{proba-diff}
\beqs
\limsup_{n\rightarrow\infty} \underline \psi_n \big((B_s)_{s\in[0,M]}\big)~~ \leq ~~ \tau  ~~ \leq ~~  \limsup_{n\rightarrow\infty} \bar \psi_n \big((B_s)_{s\in[0,M]}\big)\;,\quad \P-a.s.
\enqs
and 
\beqs
\P\Big(\limsup_{n\rightarrow\infty}\underline \psi_n\big((B_s)_{s\in[0,M]}\big)~\neq~ \limsup_{n\rightarrow\infty}\bar \psi_n^k\big((B_s)_{s\in[0,M]}\big)\Big) & = & 0
\enqs
Taking $\psi = \limsup_{n\rightarrow\infty}\bar \psi_n$ we get $\tau  =   \psi\big((B_s)_{s\in[0,M]}\big)$ $\P$-a.s.

\vspace{2mm}

We now turn to $\zeta$. Since $\zeta$ is $\Hc_\tau$-measurable, $\zeta\mathds{1}_{\tau\leq t}$ is $\Hc_t$-measurable for all $t\geq 0$. Using $\tau\leq M$ $\P$-a.s. we get
$\zeta$ is $\Hc_M$-measurable. Using Remark $32$, Chapter 2 in \cite{delmey75} as previously done, we get a Borel function $\phi$ such that
\beqs
\zeta & = & \phi\big((B_s)_{s\in[0,M+1]}\big) ~~\quad \P-a.s.
\enqs
\ep

\begin{Proposition}\label{egal-loi-meas}
Let $(\Omega^i,\Gc^i,\P^i)$, $i=1,2$, be two  compete probability spaces. Suppose that each $(\Omega^i,\Gc^i,\P^i)$ is endowed with a Brownian motion $W^i$ 
 and denote by $\F^i=(\Fc^i_t)_t$ the filtration satisfying usual conditions generated by $W^i$. 

 Fix $(\tau^i,\zeta^i)$ a couple of random variables defined on $(\Omega^i,\Gc^i,\P^i)$ for $i=1,2$ and suppose that
 \begin{itemize}
 \item  $\tau^1$ is an $\F^1$-stopping time,
\item $\zeta^1$ is $\Fc^1_{\tau^1}$-measurable  
\item $(W^2,\tau^2,\zeta^2)$ has the same law as $(W^1,\tau^1,\zeta^1)$.
\end{itemize}
Then $\tau^2$ is an $\F^2$-stopping time and $\zeta^2$ is $\Fc^2_{\tau^2}$-measurable.
\end{Proposition}
\noindent\textbf{Proof.} 
Since   $\tau^1$ is an ${\F^1}$-stopping time and ${\F^1}$ is the complete right-continuous filtration of $(W^1_s)_{s\geq 0}$, we can write from Remark $32$, Chapter 2 in \cite{delmey75}  for any $r\geq0$ and any $\eps>0$, 
\beqs
\underline \psi\big((W^1_s)_{s\in[0,r+\eps]}\big) ~\leq ~\mathds{1}_{\tau^1\leq r}~\leq~\bar \psi\big((W^1_s)_{s\in[0,r+\eps]}\big)
\enqs
and 
\beqs
\P^1\big(\underline \psi\big((W^1_s)_{s\in[0,r+\eps]}\big)~\neq~ \bar \psi\big((W^1_s)_{s\in[0,r+\eps]}\big)\big) & = & 0
\enqs
where $\underline \psi$ and $\bar \psi$ are two Borel functions. Since $(W^1,\tau^1)$ and  $(W^2, \tau^{2})$ have the same law
 we get
 \beqs
 \P^2\Big(\underline \psi\big((W^2_s)_{s\in[0,r+\eps]}\big) ~\leq ~\mathds{1}_{\tau^{2}\leq r}~\leq~\bar \psi\big((W^2_s)_{s\in[0,r+\eps]}\big)\Big) & = & 1
\enqs
and 
\beqs
 \P^2\Big(\underline \psi\big((W^2_s)_{s\in[0,r+\eps]}\big)\neq \bar \psi\big((W^2_s)_{s\in[0,r+\eps]}\big)\Big) & = & 0\;.
\enqs
Since $\F^2$ is complete this implies that $\mathds{1}_{\tau^{2}\leq r}$ is $\Fc^2_{r+\eps}$-measurable. Using the right-continuity of $\F^2$, we deduce that $\mathds{1}_{\tau^{2}\leq r}$ is $\Fc^2_{r}$-measurable and $\tau^2$ is an $\F^2$-stopping time.

  By the same argument, we get that the random variable $\zeta^{2}\mathds{1}_{\tau^{2}\leq r}$ is $\Fc^{2}_r$-measurable for all $r\geq 0$, which is equivalent to the $\Fc^{2}_{\tau^{2}}$-measurability of $\zeta^{2}$. \ep

\vspace{2mm}

We now provide two results on measurability and convergence for a sequence of processes defined on the same space but with different filtrations.

\vspace{2mm}

We fix in the sequel  a complete probability space $(\Omega,\Gc,\P)$  on which is defined a sequence of Brownian motions $(B^n)_{n\geq 0}$. For $n\geq0$, we denote by $\F^n=(\Fc^n_t)_{t\geq0}$ the complete right-continuous filtration generated by $B^n$.

\begin{Proposition}\label{prop-meas-conv} 
For $n\geq 1$, let $\tau^n$ be an $\F^n$-stopping time and  $\zeta^n$ be an $\Fc_{\tau^n}^n$-measurable random variable.
We suppose that
\begin{enumerate}[(i)]
\item $B^n$ converges to $B^0$:
\beqs
\sup_{t\in[0,T]}|B^n_t-B^0_t| & \xrightarrow[n\rightarrow\infty]{\P-a.s.}  & 0\;,
\enqs
\item the sequences $(\tau^n)_{n\geq1}$ and $(\zeta^n)_{n\geq1}$ are uniformly bounded,
\item
there exist random variables $\tau^0$ and $\zeta^0$ such that
\beqs
(\tau^n,\zeta^n) & \xrightarrow[n\rightarrow\infty]{\P-a.s.} & (\tau^0,\zeta^0)\;.
\enqs
\end{enumerate}
Then, $\tau^0$ is an $\F^0$-stopping time and $\zeta^0$ is $\Fc^0_{\tau}$-measurable.
\end{Proposition}
\ni\textbf{Proof.}
We first prove that $\tau^0$ is an $\F^0$-stopping time. Fix $t> 0$ and define for $p\geq 1$, the bounded and continuous functions $\Phi_p$ by
 \beqs
 \Phi_p(x) & = & \mathds{1}_{x\leq t-{1\over p}}+p\mathds{1}_{t-{1\over p}< x \leq t}(t-x) \;,\quad x\in \R_+
 \enqs 
 From Theorem 3.1 in \cite{BDM01}  and (iii) we get
\beqs
\E\big[\Phi_p(\tau^{n})|\Fc^{n}_t\big] & \xrightarrow[n\rightarrow\infty]{\P }& \E\big[\Phi_p({\tau^{0}})|\Fc^{0}_t\big]\;.
\enqs
 Since $\tau^n$ is an $\F^n$-stopping time we have $\E[\Phi_p(\tau^{n})|\Fc^{n}_t] = \Phi_p(\tau^{n})$. Indeed, we can write $\Phi_p=\lim_{k\infty}\Phi_p^k$ where $\Phi_p^k$ is defined by
 \beqs
\Phi_p^k(x) & = & \mathds{1}_{x\leq t-{1\over p}}+\sum_{j=1}^k{j\over kp}\mathds{1}_{t-{j\over kp}< x \leq t-{j-1\over kp}} \;,\quad x\in \R_+\;.
\enqs
 Then since $\tau^n$ is an $\F^n$ stopping time, the random variable $\Phi^k_p(\tau^n)$ is $\Fc^n_t$-measurable. Sending $k$ to infinity, we get that $\Phi_p(\tau^n)$ is $\Fc^n_t$-measurable.
 
 Since $\Phi_p$ is continuous we get from (iii)
 \beqs
\Phi_p(\tau_{n}) & \xrightarrow[n\rightarrow\infty]{ \P-a.s.}& \Phi_p(\tau^{0}).
 \enqs
 Therefore $\Phi_p(\tau^{0})=\E[\Phi_p(\tau^{0})|\Fc^{0}_t] $. Sending $p$ to infinity we get $\mathds{1}_{\tau^{0}\leq t}= \E[\mathds{1}_{\tau^{0}\leq t}|\Fc^{0}_t]$ and $\tau^0$ is a $\F^0$-stopping time since $\F^0$ is complete. 
 
 To prove that $\zeta^0$ is $\Fc^0_{\tau^0}$-measurable, we proceed in the same way and consider $\zeta^n\Phi_p(\tau^n)$ instead of $\Phi_p(\tau^n)$ for $n\geq 0$.
\ep

\vspace{4mm}
\noindent We now turn to stability of diffusions. For $n\geq 0$, we fix random functions $b_n:[0,T]\times\Omega\times\R^d\rightarrow\R^d$ and $a_n:[0,T]\times\Omega\times\R^d\rightarrow\R^{d\times d}$. We suppose that 

\vspace{2mm}

\ni\textbf{(HA)}
\begin{enumerate}[(i)]
\item For each $n\geq 0$, $b_n$ and $a_n$ are $\F^n$-progressive$\otimes\Bc(\R^d)$-measurable,
\item there exists $\delta>0$ such that 
\beqs
\E\Big[\int_0^T\Big(|b^n(t,0)|^{2+\delta}+|a^n(t,0)|^{2+\delta}\Big)dt \Big]& < & +\infty\;,
~n\geq 0\;,
\enqs
\item there exists a constant $L$ such that
\beqs
|b^n(t,x)-b^n(t,x')|+|a^n(t,x)-a^n(t,x')| & \leq & L|x-x'|\;,\quad x,x'\in\R^d\;,~n\geq 0\;.
\enqs
 \end{enumerate}
 Then, for a given deterministic initial condition $X_0$, we can define for each $n\geq 0$, the solution $X^n$ to the SDE
 \beqs
 X^n_t & = & X_0+\int_0^tb^n(s,X^n_s)ds+\int_0^ta^n(s,X^n_s)dB^n_s \qquad t\geq0\;.
 \enqs
\begin{Proposition}\label{AnConvDiff} Suppose that 
 \beq\label{convMBSB}
\sup_{t\in[0,T]}|B^n_t-B^0_t| & \xrightarrow[n\rightarrow\infty]{\P-a.s.}  & 0\;,
\enq
and
\beq\label{estimbnb}
\E\Big[\int_0^T\big|a^n(s,x)-a^0(s,x)\big|^2ds\Big]+\E\Big[\int_0^T|b^n(s,x)-b^0(s,x)|^2ds\Big] & \xrightarrow[n\rightarrow+\infty]{} & 0\;,\qquad 
\enq
for all $x\in\R^d$. Then, under \textbf{(HA)}, we have
\beq
\E\Big[ \sup_{t\in[0,T]}\big|X^n_t-X^0_t\big|^{2} \Big] & \xrightarrow[n\rightarrow\infty]{} & 0\;.\label{convXnX}
\enq
\end{Proposition}
To prove this result we cannot use classical estimates on diffusions processes since the driving Brownian motion evolves with $n$. In particular the stochastic integrals $\int a^ndB^0$ are not defined. We therefore need to use approximations by step processes as done in the construction of the It\^o integral.

\vspace{2mm}

\noindent\textbf{Proof.} 
We proceed in two steps.

\vspace{2mm}

\noindent \textbf{Step 1.} We first consider the case where the $b^n$ and $a^n$ do not depend on the variable $x$. 
For $p\geq1$, Let $H^p$ be an $\F$-adapted piecewise constant process 
of the form
\beqs
H^p_t & = & \sum_{k=0}^{N_p} \tilde H^p_k\mathds{1}_{[t^p_k,t^p_{k+1})}(t)\;,\quad t\in [0,T] 
\enqs
where $\tilde H^p_k\in \mathbf{L}^{2+\delta}(\Omega,\Fc_{t^p_k},\P)$ for $0\leq k\leq N_p$, such that
\beq\label{estimHpa}
\E\Big[\int_0^T |H^p_s-a_s|^2ds\Big] & \leq & {1\over p}\;.
\enq
We then have 
\beq\label{anmoinsa}
\E\Big[\Big|\int_0^Ta^ndB^n-\int_0^TadB^0\Big|^2\Big] & \leq & 2\Big(E\Big[\Big|\int_0^Ta^ndB^n-\int_0^TH^pdB^0\Big|^2\Big] +{1\over p} \Big)\;.
\enq
We then define the process $H^{p,n}$ by
\beqs
H^{p,n}_t & = & \sum_{k=0}^{N_p} \E\Big[\tilde H^p_k\big|\Fc^n_{t^p_k}\Big]\mathds{1}_{[t^p_k,t^p_{k+1})}(t)\;,\quad t\in [0,T] \;.
\enqs
We can write the following decomposition
\beq
\E\Big[\Big|\int_0^Ta^ndB^n-\int_0^TH^pdB^0\Big|^2\Big] & \leq & 2\Big(E\Big[\Big|\int_0^Ta^ndB^n-\int_0^TH^{p,n}dB^n\Big|^2\Big] \nonumber  \\
 & & +E\Big[\Big|\int_0^TH^{p,n}dB^n-\int_0^TH^pdB^0\Big|^2\Big] \Big)\;.\label{anHp}
\enq
From  \reff{convMBSB}, we can apply Proposition 2 in \cite{CMS} and we get
\beq\label{convHnpk}
\E\Big[\tilde H^p_k\big|\Fc^n_{t^p_k}\Big] & \xrightarrow[n\rightarrow+\infty]{\P} & \tilde H^{p}_k \;,\quad 0\leq k\leq N_p\;. 
\enq 
In particular we get from \reff{convMBSB} and \reff{convHnpk}
\beq\label{HnpBnHpB}
\E\Big[\Big|\int_0^TH^{p,n}dB^n-\int_0^TH^pdB^0\Big|^2\Big] & \xrightarrow[n\rightarrow+\infty]{} & 0\;.
\enq
Moreover, from It\^o Isometry and \reff{estimHpa} we have
\beq
E\Big[\Big|\int_0^Ta^ndB^n-\int_0^TH^{p,n}dB^n\Big|^2\Big]  & = & E\Big[\int_0^T\big|a^n_s-H^{p,n}_s\big|^2ds\Big] \nonumber \\
& \leq  & 
3\Big(E\Big[\int_0^T\big|a^n_s-a^0_s\big|^2ds\Big]+{1\over p}\label{decompanBnHpnBn}
 \\& & \nonumber  \hspace{6mm}
 +E\Big[\int_0^T\big|H^{p}_s-H^{p,n}_s\big|^2ds\Big] \Big)\;.
\enq
Then using \reff{convHnpk}, we also get
\beq\label{estimHpHpn}
\E\Big[\int_0^T\big|H^{p}_s-H^{p,n}_s\big|^2ds\Big] & \xrightarrow[n\rightarrow+\infty]{} & 0\;.
\enq
Therefore, we get from \reff{estimbnb}, \reff{decompanBnHpnBn} and \reff{estimHpHpn}
\beqs
\limsup_{n\rightarrow\infty}\E\Big[\Big|\int_0^Ta^ndB^n-\int_0^TH^{p,n}dB^n\Big|^2\Big]  & \leq & {1\over p}\;.
\enqs
From this last inequality, \reff{anmoinsa}, \reff{anHp} and \reff{HnpBnHpB} we get
\beqs
\limsup_{n\rightarrow\infty} 
\E\Big[\Big|\int_0^Ta^ndB^n-\int_0^Ta^0dB^0\Big|^2\Big] & \leq & {4\over p}\;,  \quad p\geq1\;.
\enqs
Therefore, we get
\beqs
\lim_{n\rightarrow\infty} 
\E\Big[\Big|\int_0^Ta^ndB^n-\int_0^Ta^0dB^0\Big|^2\Big] & = & 0\;.
\enqs
From Theorem 3.1 in \cite{BDM01}, we deduce that
\beqs
\lim_{n\rightarrow\infty} 
\E\Big[\sup_{t\in[0,T]}\Big|\int_0^ta^ndB^n-\int_0^ta^0dB^0\Big|^{2}\Big] & = & 0\;.
\enqs
 From this last equality and \reff{estimbnb}, we get \reff{convXnX}.

\vspace{2mm}

\ni \textbf{Step 2.} We now consider the general case. For $n\geq 0$, we denote by $(X^{n,p})_{p\geq0}$ the sequence of processes defined by
\beqs
X^{n,0}_t  & = & X_0\;,\quad t\geq 0\;,
\enqs 
and 
\beqs
X^{n,p+1}_t  & = & X_0+\int_0^tb^n(s,X^{n,p}_s)ds+\int_0^ta^n(s,X^{n,p}_s)dB^n_s\;,\quad t\geq0\;,
\enqs 
for $p\geq 0$. From \textbf{(HA)} (ii) and since $X_0$ is deterministic, we get by induction on $p$ that 
\beqs
\E\Big[\sup_{t\in[0,T]} |X^{n,p}_t|^{2+\delta}\Big] & < & \infty
\enqs
for all $n,p\geq 1$.  Still using an induction we get from Step 1 that
 \beq
\E\Big[ \sup_{t\in[0,T]}\big|X^{n,p}_t-X^{0,p}_t\big|^2 \Big] & \xrightarrow[n\rightarrow\infty]{} & 0\;\label{convXnpX0p}
\enq
for all $p\geq 0$.
From argument on diffusion processes, we have (see e.g. the proof of Theorem 2.9 of Chapter 5 in \cite{KS91})
\beqs
\sup_{n\geq 0} \E\Big[\sup_{t\in[0,T]}|X^{n,p}_t-X^{n}_t|^2\Big] & \leq  & \psi(p)
\enqs
where $\psi(p)\rightarrow0$ as $p\rightarrow+\infty$. We then get
\beqs
\limsup_{n\rightarrow+\infty}\E\Big[ \sup_{t\in[0,T]}\big|X^n_t-X^0_t\big|^2 \Big]  & \leq & 2\psi(p)+\lim_{n\rightarrow+\infty}\E\Big[ \sup_{t\in[0,T]}\big|X^{n,p}_t-X^{0,p}_t\big|^2 \Big]~~\leq~~2\psi(p)\;.
\enqs
Sending $p$ to $\infty$, we get the result.
\ep
\subsection{Proofs of Lemmata \ref{Propfngn} and \ref{prop-bound-v}}\label{appendice preuve lemme}

\ni\textbf{Proof of Lemma \ref{Propfngn}.}   Fix $n\geq1$, $R>0$ and $i\in\Ic$. From the definition of $f_n$ we have 
\beqs
|f_{n}(x,i)-f_{n}(x',i)| & \leq & n\big|\Theta_{n}(x)-\Theta_{n}(x')\big|+ \big|f(x,i)-f(x',i)\big|\,,
\enqs
for all $x,x'\in\R^d$ and $i\in\Ic$.
Since $d(.,\Dc)$ is Lipschitz continuous, we get from the definition of $\Theta_n$ and \textbf{(H2)} (i) the existence of a constant $L_{R,n}$ such that
\beqs
|f_{n}(x,i)-f_{n}(x',i)|
& \leq & L_{R,n}|x-x'|\,,
\enqs
for all $x,x'\in\R^d$. 

We turn to the grow property. From the definition of $f_n$ we have 
\beqs
|f_{n}(x,i)| & \leq & n\big|\Theta_{n}(x)\big|+ \big|f(x,i)\big|\,,
\enqs
for all $x\in\R^d$ and $i\in\Ic$.
Since $d(.,\Dc)$ is Lipschitz continuous, it has a linear growth and we get from the definition of $\Theta_n$ and \textbf{(H2)} (ii) that there exists a constant $C_n$ such that
\beqs
|f_{n}(x,i)| & \leq & C_n\big(1+ |x|^q\big)\,,
\enqs
The proof is the same for the function $g_{n}$.
\ep

\vspace{2mm}

\ni\textbf{Proof of Lemma  \ref{prop-bound-v}.} Fix $n\geq 1$ and $(t,x,i)\in[0,T]\times\Dc\times\Ic$. Using the definition of $f_n$ and $g_n$ we have
\beq\label{malJnJ1}
J_n(t,x,\alpha) & \leq & J_1(t,x,\alpha)
\enq
for any $\alpha\in\Ac_{t,i}$.
From \reff{XS2} and \textbf{(H2)} there exists a constant $C$ such that
\beqs
J_1(t,x,\alpha) & \leq & C\big(1+|x|^q\big)
\enqs
for any $\alpha\in\Ac_{t,i}$. From \reff{malJnJ1} and the definition of $v_n(t,x,i)$, we get $\reff{maj-lin-v}$.
\ep

\vspace{13mm}

\begin{small}

\end{small}

\end{document}